\RequirePackage{ifpdf}
\ifpdf 
\documentclass[pdftex]{sigma}
\else
\documentclass{sigma}
\fi

\usepackage[all]{xy}
\numberwithin{equation}{section}

\def\Cal{\mathcal}

\def\b1{\text{\bf 1}}

\def\BC{{\Bbb C}}
\def\BP{{\Bbb P}}
\def\BV{{\Bbb V}}

\def\BP{{\Bbb P}}

\def\BZ{{\Bbb Z}}

\def\CD{{\Cal D}}

\def\CH{{\Cal H}}
\def\CL{{\Cal L}}

\def\CN{{\Cal N}}
\def\CO{{\Cal O}}

\def\CT{{\Cal T}}
\def\CV{{\Cal V}}
\def\CU{{\Cal U}}
\def\CW{{\Cal W}}

\def\Der{\mathrm{Der}}

\def\ft{{\frak t}}

\def\fg{{\frak g}}
\def\fa{{\frak a}}
\def\fn{{\frak n}}

\def\ft{{\frak t}}
\def\fq{{\frak q}}

\def\Ker{\text{Ker}}

\def\Vir{{\Cal V}{\rm ir}}
\def\Vert{{\Cal V}{\rm ert}}

\def\iso{ \buildrel\sim\over\longrightarrow }

\begin{document}

\allowdisplaybreaks

\renewcommand{\thefootnote}{$\star$}

\renewcommand{\PaperNumber}{086}

\FirstPageHeading

\ShortArticleName{Vertex Algebroids over Veronese Rings}

\ArticleName{Vertex Algebroids over Veronese Rings\footnote{This paper is a
contribution to the Special Issue on Kac--Moody Algebras and Applications. The
full collection is available at
\href{http://www.emis.de/journals/SIGMA/Kac-Moody_algebras.html}{http://www.emis.de/journals/SIGMA/Kac-Moody{\_}algebras.html}}}

\Author{Fyodor MALIKOV}

\AuthorNameForHeading{F.~Malikov}

\Address{Department of Mathematics, University of Southern California, Los Angeles, CA 90089, USA}
\Email{\href{mailto:fmalikov@usc.edu}{fmalikov@usc.edu}}

\ArticleDates{Received July 28, 2008, in f\/inal form December 07,
2008; Published online December 13, 2008}

\Abstract{We f\/ind a canonical quantization of Courant algebroids
over Veronese rings. Part of our approach allows a semi-inf\/inite
cohomology interpretation, and the latter can be used to def\/ine
sheaves of chiral dif\/ferential operators on some homogeneous spaces
including the space of pure spinors punctured at a point.}


\Classification{14Fxx, 81R10; 17B69}

\section{Introduction}

Attached to a commutative associative algebra $A$  are the Lie
algebra of its
 derivations, $\Der(A)$, and the module of K\"ahler dif\/ferentials,
 $\Omega(A)$. The  identities that are
 satisf\/ied by the classic dif\/ferential geometry operations, such as
 the Lie bracket, the Lie derivative, the de Rham dif\/ferential,
 etc., can be summarized by saying that the $A$-module $\Der(A)\oplus\Omega(A)$
 is a {\it Courant algebroid}, \cite{C,LWX}. For  reasons that will become
 apparent later, we will use the notation
 $\CV^{\rm poiss}(A)=\CV^{\rm poiss}(A)_0\oplus\CV^{\rm poiss}(A)_1$, where
 $\CV^{\rm poiss}(A)_0=A$ and $\CV^{\rm poiss}(A)_1=\Der(A)\oplus\Omega(A)$.

 This example can be enriched in two dif\/ferent ways. First, it can
 be quantized. Attached to~$A$ in~\cite{GMSI}  is the notion of a
 {\it vertex algebroid}, $\CV(A)$. This notion is a result of
 axiomatizing the   structure that is
 induced on conformal weight 0 and 1 components of a graded vertex
 algebra. One has $\CV(A)=A\oplus\CV(A)_1$ for some $\CV(A)_1$,
 which f\/its in the exact sequence
 \[
 0\rightarrow\Omega(A)\rightarrow\CV(A)_1\rightarrow
 \Der(A)\rightarrow 0.
 \]
 Therefore, $\CV(A)_1$ is f\/iltered and the corresponding graded
 object is ${\rm Gr}\, \CV(A)_1=\CV^{\rm poiss}(A)_1$. This strongly resembles the
 Poincar\'e--Birkhof\/f--Witt f\/iltration, and it is indeed true that the
 notion of a Courant algebroid is a quasiclassical limit of that of
 a vertex algebroid; this important observation is due to Bressler~\cite{Bre}, but the fact that the Courant bracket belongs in the inf\/inite
 dimensional world had been discovered by Dorfman much earlier,~\cite{Dor}. The  relation of this notion to various string theory
 models has been elucidated in~\cite{Mal}.

 Unlike its quasiclassical counterpart, a vertex algebroid
 may not exist, and if it exists, it may not be unique.
 If ${\rm Spec}(A)$ is smooth, then \cite{GMSI}  the obstruction to existence is the class
 ${\rm ch}_2({\rm Spec}(A))$, and if ${\rm ch}_2({\rm Spec}(A))=0$, then
  the isomorphism classes are
 parameterized by the hypercohomology group $H^1({\rm Spec}(A),
 \Omega_{A}^{2}\stackrel{d_{DR}}{\rightarrow}\Omega_{A}^{3,{\rm cl}})$.
 But what if ${\rm Spec}(A)$ is not smooth?

This question suggests the second way to enrich, which is to note
that $\Der(A)$ and $\Omega(A)$ tell the whole story only if ${\rm Spec}(A)$
is smooth; if not, then the higher algebras of derivations and
modules of 1-forms must arise. This was made precise by Hinich,
\cite{Hin}, who def\/ined $\Der(A)^{\bullet}$ as a kind of a derived
functor of the functor $\Der$ by applying the latter to a polynomial
dif\/ferential graded algebra resolution $R\rightarrow A$. This gives
rise to a graded Courant algebroid functor
$A\mapsto\CV^{\rm poiss}(A)^{\bullet}$.

The problem of combining the two, that is to say, f\/inding a
quantization, $\CV(A)^{\bullet}$, of $\CV^{\rm poiss}(A)^{\bullet}$
appears to be intellectually attractive and important for
applications. If $A$ is a complete intersection, then the resolving
algebra $R$ can be chosen to be a super-polynomial ring on f\/initely
many generators, the corresponding resolution $R\rightarrow A$ being
none other than the standard Koszul complex, and the quantization, a
dif\/ferential graded vertex algebroid $\CV(R)$, is immediate; this
observation has been used in a number of physics and mathematics
papers. If, however, $A$ is not a complete intersection, then any
resolving algebra $R$ is inf\/initely generated in which case def\/ining
a vertex algebroid $\CV(R)$ becomes problematic because of various
divergencies. A~re\-gu\-la\-rization procedure for some of these
divergencies was suggested in \cite{BN} and elaborated on in~\cite{GS}.

Here is what we do in the present paper. Let $V_N$ be the
$(N+1)$-dimensional irreducible $sl_2$-module, $O_N\subset\BP(V_N)$
the highest weight vector orbit, and $A_N$ the corresponding
homogeneous coordinate ring. All of this is a representation
theorist's way of saying that $\BP^1\iso O_N\subset\BP^N$ is a
Veronese curve, and $A_N$ is a Veronese ring.

$A_N$ is a quadratic algebra, in fact it is Koszul \cite{BF,IM,
Bez}, but it is not a complete intersection. The main result of the
paper, Theorem~5.1.1, 
asserts that
$\CV^{\rm poiss}(A_N)^{\bullet}$ admits a unique quantization. It is no
surprise that this quantization, $\CV(A_N)$, contains a vertex
algebroid, $\CV(sl_2)_{k}$,  attached to $\widehat{sl}_2$ with some
central charge  $k$. What is more important is that the vertex
algebroid attached to $\widehat{gl}_2$ enters the fray. The latter,
$\CV(gl_2)_{k_1 k_2}$, depends in general on two central charges,
$k_1$, $k_2$, and we f\/ind that the quantization conditions imply,
f\/irst, that $k_1+k_2=-2$ and, second, that $k_1=-N-2$.

Theorem 5.1.1
and its proof appear in Section~\ref{sect.4}, and it is for the sake of this
section that the paper was written. Section~\ref{sect.3} is to a large extent an exposition of
Hinich's result (see also \cite{Behr}) with some extensions (Sections~4.3, 
4.4, 
4.5) 
that are needed in Section~\ref{sect.4}. Sections~\ref{sect.1} and
\ref{sect.2} are an attempt, perhaps futile, to make the paper self-contained~-- except
Sections~3.7.3, 
3.7.4, 
where sheaves of vertex algebroids over
$\BC^2\setminus 0$ are classif\/ied. The classif\/ication obtained is instrumental in proving
Theorem~5.1.1; 
in particular, the vertex algebroid $\CV(gl_2)_{k_1,k_2}$ with the
compatibility condition $k_1+k_2=-2$ makes appearance in Section~3.7.4. 

An obvious generalization of $A_N$ is provided by the homogeneous
coordinate ring of a higher dimensional Veronese embedding
$\BP(\BC^n)\rightarrow \BP(S^N(\BC^n))$. We show (Theorem~5.3.1) 
that if $n>2$ and $N>1$, then no quantization
exists.

Much of the above carries over to an arbitrary simple $\fg$, where
$A_N$ is replaced with the homogeneous coordinate ring of the
highest weigh vector orbit in the projectivization of a simple
$\fg$-module. For example, $\BC^2\setminus 0$ becomes the
Bernstein--Gelfand--Gelfand base af\/f\/ine space, $G/N$.  Constructed in~\cite{GMSII} is the 1-parameter family of sheaves of vertex
algebroids $H^{\infty/2}(L\fn,\CV_{G,k})$ over~$G/N$, where
$\CV_{G,k}$, $k\in\BC$, is a family of vertex algebroids over~$G$,
\cite{AG, F1, FP, GMSII}. There is little doubt that the family
$H^{\infty/2}(L\fn,\CV_{G,k})$, $k\in\BC$, is universal in that it
classif\/ies vertex algebroids over $G/N$ equipped with
$\CV(\fg)_k$-structure. This is a higher rank analogue of the
classif\/ication obtained in Section~3.7.3 
and alluded to
above. Note that just as $G/N$ is a $G\times T$-space, the maximal
torus acting on the right, so there is a diagram of embeddings
\[
\CV(\fg)_{k_1}\hookrightarrow
H^{\infty/2}(L\fn,\CV_{G,k})\hookleftarrow \CV(\ft)_{k_2}\qquad \text{with} \quad k_1+k_2=-h\check{}\,.
\]
Therefore, $\CV(\fg)_{k}\oplus \CV(\ft)_{-k-h \check{ }}$ is a
higher rank analogue of $\CV(gl_2)_{k,-k-2}$; here $\fg$, $\fn$ and
$\ft$ are the Lie algebras of the Lie groups $G$, $N$, and $T$ resp.
We elaborate on these remarks in Section~\ref{sect.5}, where we use the technique
of semi-inf\/inite cohomology to compute CDO-s on some homogeneous
spaces including the spaces of pure spinors punctured at a point. In
the latter case, this gives an approach alternative to that of the
original result by Nekrasov \cite{N}.

Some aspects of the $sl_2$-case, however, are not that easy to
generalize. As they say, we hope to return to this subject in a
separate paper.

We would like to conclude by saying that a major source of
inspiration was provided to us by the work of Berkovits and Nekrasov
\cite{Ber,BN}, where similar problems are analyzed in the case of
the spinor representation of the spinor group.

\section{Vertex algebras}\label{sect.1}

\paragraph{2.1. Conventions.} 
Underlying all the
constructions in this paper will be the category of $\BZ$-graded
vector superspaces and grading preserving linear maps over $\BC$.
This grading will be called (and should be thought of as) the {\it
homological degree} grading. More often, though, the attribute
`graded' will be skipped. Thus the phrase `let $V$ be a vector
space' will mean that $V=\oplus_{n\in\BZ}V^{n}$,
$V^{\rm even}=\oplus_{n\in\BZ}V^{2n}$,
$V^{\rm odd}=\oplus_{n\in\BZ}V^{2n+1}$. Likewise, the pref\/ix `super-'
will be usually omitted so that commutative will mean
super-commutative, algebra super-algebra, bracket super-bracket:
$[a,b]=ab-(-1)^{ab}ba$.

If $V$ and $W$ are vector spaces, then $V\otimes W$ is also a vector
space with  homological degree grading def\/ined in the standard way
so that $(V\otimes W)^n=\oplus_{i\in\BZ}V^i\otimes W^{n-i}$. Various
bilinear operations (`multiplications') to be used below will be
morphisms of graded vector spaces $V\otimes W\rightarrow U$.

Along with the homological degree grading, the {\it grading by
conformal weight} will play a~prominent role. The latter will be
indicated by a subindex; thus, for example, the phrase `a~graded (by
conformal weight) vertex algebra' will mean, in particular, a vector
space $V$ with a~direct sum decomposition $V=\oplus_{n\geq 0}V_n$
valid in the category of graded vector spaces.

Most of the def\/initions and constructions in this and the following
section are well known, and their graded versions are always
straightforward. We recommend \cite{Kac} and \cite{FBZ} as an
excellent introduction to vertex things  and a guide to further
reading.

 \paragraph{Def\/inition 2.2.} 
  A vertex algebra is a collection $(V,\b1, T,
_{(n)}, n\in\BZ)$, where $V$
 is a vector space,  $\b1\in V$  is a distinguished element known as the {\it vacuum
 vector}, $T:V\rightarrow V$ is a linear operator known as the {\it
 translation operator}, each $_{(n)}$ is a multiplication
 \[
 _{(n)}: \ \ V\otimes V\rightarrow V,\quad  a\otimes b\mapsto a_{(n)}b \quad \text{s.t.}
 \quad a_{(n)}b=0\quad \text{if} \quad
 n\gg 0,
 \]
 that is subject to the following axioms:

(1) (vacuum)
\begin{gather}
\b1_{(n)}=\left\{\begin{array}{ll}
{\rm Id}_{V}&\text{if} \ \ n=-1, \\
0&\text{otherwise},
\end{array}
\right.  \qquad  a_{(-1)}\b1=a,\quad \forall\;  a\in V; \label{(1.1a)}
\end{gather}

(2) (translation invariance)
\begin{equation}
[T,a_{(n)}]b=(Ta)_{(n)}b=-na_{(n-1)}b,\quad \forall \; a,b\in V, \quad
n\in\BZ;\label{(1.1b)}
\end{equation}

(3) (skew-symmetry)
\begin{equation}
a_{(n)}b=(-1)^{ab}\sum_{j\geq 0}(-1)^{n+1+j}
\frac{1}{j!}T^{j}(b_{(n+j)}a),\quad \forall \; a,b\in V, \quad
n\in\BZ;\label{(1.1c)}
\end{equation}

(4) (Jacobi identity)
\begin{equation}
[a_{(m)},b_{(n)}]c=\sum_{j\geq
0}\binom{m}{j}(a_{(j)}b)_{(m+n-j)}c,\quad \forall\;  a,b,c\in V,\quad
m,n\in\BZ; \label{(1.1d)}
\end{equation}

(5) (quasi-associativity or normal ordering)
\begin{gather}
(a_{(-1)}b)_{(n)}c= \sum_{j\geq
0}a_{(-1-j)}b_{(n+j)}c+(-1)^{ab}\sum_{j> 0}b_{(n-j)}a_{(-1+j)}c,\quad
\forall \; a,b,c\in V,\quad  n\in\BZ.\label{(1.1e)}
\end{gather}

The collection of axioms we used in Def\/inition~2.2 
 is a
little redundant but makes the exposition a little more transparent.
It emphasizes the fact that the notion of a vertex algebra is a
mixture of (appropriate analogues of) that of an associative algebra
and a Lie algebra. Extracting the Lie part of the def\/inition one
arrives at the notion of a vertex Lie algebra.

\paragraph{Def\/inition 2.3.} 
A {\it vertex Lie algebra} is a
collection $(V, T,{}_{(n)},n\in\BZ_{+})$, where $V$ is a vector
space, $T:V\rightarrow V$ is a linear operator, each $_{(n)}$ is a
multiplication
\begin{equation}
 _{(n)}: \ \ V\otimes V\rightarrow V,\quad  a\otimes b\mapsto a_{(n)}b \quad \text{s.t.} \quad a_{(n)}b=0
 \quad \text{if} \quad n\gg 0
\label{(1.2)}
\end{equation}
 that is subject to the following axioms:

(1) translation invariance, that is, (\ref{(1.1b)}) for $n\geq 0$;

(2) skew-symmetry, that is, (\ref{(1.1c)}) for $n\geq 0$;

(3) Jacobi identity, that is, (\ref{(1.1d)}) for $n\geq 0$.

\medskip

There is an obvious forgetful functor
\begin{equation}
\Phi:  \ \ \{\text{Vertex algebras}\}\rightarrow \{\text{Vertex Lie
algebras}\}.\label{(1.3a)}
\end{equation}
Its left adjoint functor (the vertex enveloping algebra functor)
\begin{equation}
\CU: \ \ \{\text{Vertex Lie algebras}\}\rightarrow \{\text{Vertex
algebras}\}\label{(1.3b)}
\end{equation}
is well known to exist, see~\cite{Pr};  it also appears in
\cite{Kac} as `the vertex algebra attached to a formal distribution
Lie superalgebra'.

Note a canonical map
\begin{equation}
\iota: \ \ \CL\rightarrow \Phi(\CU\CL)\label{(1.4)}
\end{equation}
that is the image of ${\rm Id}\in {\rm Hom}(\CU\CL,\CU\CL)$ under the
identif\/ication ${\rm Hom}(\CU\CL,\CU\CL) \!\iso\! {\rm Hom}(\CL,\Phi\CU\CL)$.

\paragraph{Example 2.4.} 
Let $\CL'$ be a free
$\BC[T]$-module on one generator $L$ and let
$\CL(\Vir)_{c}=\CL'\oplus\BC$, where $\BC$ is considered as a
trivial $\BC[T]$-module. $\CL(\Vir)_{c}$ carries a unique vertex Lie
algebra structure such that
\[
L_{(0)}L=T(L),\quad L_{(1)}L=2L,\quad L_{(2)}L=0,\quad
L_{(3)}L=\tfrac{1}{2}c,\quad L_{(n)}L=0 \quad \text{if} \ \ n>3.
\]
 Upon quotienting out by the relation $\b1=1$, the vertex enveloping algebra $\CU\CL(\Vir)_{c}$ becomes
  the vacuum
representation of the Virasoro algebra of central charge $c$.

\paragraph{Example 2.5.} 
Let $\fg$ be a  Lie
algebra with an invariant bilinear form $(\cdot,\cdot)$. Let
\[
\CL(\fg)_{k}=\BC[T]\otimes\fg\oplus\BC.
\]
This space carries an obvious action of $T$, where again we consider
$\BC$ as a trivial $\BC[T]$-module, and a unique vertex Lie algebra
structure such that
\begin{equation}
(1\otimes a)_{(0)}(1\otimes b)=1\otimes [a,b],\qquad (1\otimes
a)_{(1)}(1\otimes b)=k(a,b). \label{(1.5a)}
\end{equation}
Upon quotienting out by the relation $\b1=1$, the vertex enveloping
algebra $\CU\CL(\fg)_{k}$ becomes the vacuum representation of the
corresponding af\/f\/ine Lie algebra of central charge $k$.

If $\fg$ is chosen to be $gl_{N}=sl_{N}\oplus\BC \cdot I$, then this
construction has the following version: we let
$(a,b)=\text{tr}(a\cdot b)$,
$\CL(gl_N)_{k_1,k_2}=\CL(sl_N)_{k_1}\oplus\BC[t]\otimes \BC\cdot I$
and extend (\ref{(1.5a)}) by
\begin{equation}
(1\otimes I)_{(1)}(1\otimes I)=k_{2}N,\qquad (1\otimes I)_{(0)}(1\otimes
I)=(1\otimes I)_{(n)}(1\otimes sl_N)=0\quad \forall\;  n.\label{(1.5b)}
\end{equation}
In order to handle the case of the trivial bilinear form $(\cdot,\cdot)$, or
more generally the case where $(\cdot,\cdot)$ is not unique even up to
proportionality, we will change the notation and denote by
$\CL(\fg)_{(\cdot,\cdot)}$ the vertex Lie algebra which is precisely
$\CL(\fg)_{k}$ except that the last of conditions (\ref{(1.5a)}) is
replaced with
\begin{equation}
 (1\otimes
a)_{(1)}(1\otimes b)=(a,b) \label{(1.5c)}
\end{equation}
for some $(\cdot,\cdot)$.

\medskip

A passage to the quasiclassical limit is a gentler way to blend the
Lie and commutative/as\-so\-cia\-tive algebra parts of the structure.

 \paragraph{Def\/inition 2.6.} 
 A {\it vertex Poisson algebra} is a
collection $(V,\b1\in V, T,{}_{(n)},n\geq -1)$, where $V$ is a
vector space, $T:V\rightarrow V$ is a linear operator, each ${}_{(n)}$
is a multiplication
\[
{}_{(n)}: \ \ V\otimes V\rightarrow V,\quad  a\otimes b\mapsto a_{(n)}b \quad \text{s.t.} \quad
a_{(n)}b=0 \quad \text{if} \quad n\gg0
\]
 that is subject to the following axioms:

 (1) the triple $(V,\b1, T,{}_{(-1)})$ is a unital commutative associative algebra with
 derivation;

 (2) the collection $(V,\b1\in V, T,{}_{(n)},n\geq 0)$ is a vertex Lie algebra;

 (3) each multiplication ${}_{(n)}$, $n\geq 0$, is a derivation of ${}_{(-1)}$.

 \medskip

 Vertex Poisson algebras are to vertex algebras what Poisson algebras are to noncommutative
 algebras. The following construction (cf.~\cite{Li}) is meant to illustrate this point.

 \paragraph{2.7.~Vertex algebras with f\/iltration.} 
 Suppose  a vertex algebra $V$ carries an exhaustive
increasing f\/iltration by vector
 spaces
\begin{equation}
\BC\b1\in F^{0}V\subset F^{1}V\subset\cdots\subset
F^{p}V\subset\cdots,\qquad \cup_{p\geq 0}F^{n}V=V \label{(1.6)}
\end{equation}
which satisf\/ies $ (F^{p}V)_{(n)}F^{q}V\subset F^{p+q}V$ and
$T(F^{p}V)=F^{p}V$ for all $p,q\in\BZ$ so that $
(F^{p}V)_{(n)}F^{q}V$ $\subset F^{p+q-1}V$ if $n\geq 0$. By focusing on
symbols one discovers that the vertex algebra structure on~$V$
def\/ines the following on the corresponding graded object
${\rm Gr}\, V=\oplus_{p}F^{p}V/F^{p-1}V$:
\begin{gather*}
\b1^{\rm gr}=\b1\in F^{0}V;
\\
T^{\rm gr}: \ \ F^{p}V/F^{p-1}V\rightarrow F^{p}V/F^{p-1}V,\qquad
T(\bar{a})=T(a)\mod F^{p-1}V;
\\
{}_{(-1)^{\rm gr}}: \ \ (F^{p}V/F^{p-1}V)\otimes(F^{q}V/F^{q-1}V) \rightarrow
F^{p+q}V/F^{p+q-1}V,\\
\phantom{{}_{(-1)^{\rm gr}}: \ \ }{} \ \bar{a}\otimes\bar{b} \mapsto \bar{a}_{(-1)^{gr}}\bar{b}=a_{(-1)}b\mod F^{p+q-1}V;
\\
{}_{(n)^{\rm gr}}: \ \ (F^{p}V/F^{p-1}V)\otimes(F^{q}V/F^{q-1}V) \rightarrow
F^{p+q-1}V/F^{p+q-2}V, \\
\phantom{{}_{(n)^{\rm gr}}: \ \ }{} \ \bar{a}\otimes\bar{b} \mapsto \bar{a}_{(n)^{gr}}\bar{b}=a_{(n)}b\mod F^{p+q-2}V\quad \text{if} \quad n\geq 0.
\end{gather*}
It is then immediate to check that $({\rm Gr}\, V,\b1^{\rm gr}, T^{\rm gr},{}_{(n)^{\rm gr}},n\geq -1)$ is a vertex Poisson algebra. For example,
commutativity of the product $_{(-1)}$ follows from the $n=-1$ case
of (\ref{(1.1c)}) and associativity from the $n=-1$ case of
(\ref{(1.1e)}).

The vertex algebras reviewed in Examples~2.4, 2.5 possess a f\/iltration with the indicated properties~--
as does any vertex enveloping algebra:
\begin{itemize}\itemsep=0pt
\item in  the case of $\CU\CL(\Vir)_c$ the f\/iltration is
determined by assigning degree one to $\iota(\CL(\Vir)_c)$, see
\eqref{(1.4)};

\item in  the case of $\CU\CL(\fg)_{k}$ the f\/iltration is
determined by assigning degree one to $\iota(\CL(\fg)_{k})$.
\end{itemize}

Denote  thus def\/ined vertex Poisson algebras as follows:
\begin{gather}
\CU^{\rm poiss}\CL(\Vir)={\rm Gr}\, \CU\CL(\Vir)_c,\qquad \CU^{\rm poiss}\CL(\fg)={\rm Gr}\,
\CU\CL(\fg)_{k},\nonumber\\
\CU^{\rm poiss}\CL(gl_N)={\rm Gr}\, \CU\CL(gl_N)_{k_,k_2}.
 \label{(1.7)}
\end{gather}

\section{Courant and vertex algebroids}\label{sect.2}

\paragraph{Def\/inition 3.1.} 
A vertex (vertex
Poisson) algebra $V$ is called {\it graded} if
\begin{gather}
V=\bigoplus_{n=-\infty}^{+\infty}V_{n} \quad \text{so that}\nonumber
\\
V_{n}=0,\quad {if} \ \  n<0,\ \  \b1\in V_{0},\qquad  T(V_{m})\subset
V_{m+1}, \quad \text{and} \quad (V_{m})_{(j)}V_{n}\subset
V_{m+n-j-1}.\label{(2.1)}
\end{gather}

Such grading is usually referred to as {\it conformal}, $V_n$ is
called the conformal weight $n$ component, and $v\in V_n$ is said to
have conformal weight $n$; the conformal weight of $v\in V_n$ is
usually denoted by~$\Delta(v)$.

All the vertex algebras we have seen are graded:
$\CU\CL(\Vir)_{c}$ is graded by letting $V_0=\BC$, $V_1=\{0\}$,
$V_2=\iota(L)$, $\CU\CL(\fg)_{k}$ by letting $V_0=\BC$,
$V_1=\iota(\fg)$, $\CU\CL(gl_{N})_{k_1,k_2}$ by letting
$V_1=\iota(gl_N)$

A graded vertex (vertex Poisson) algebra structure on $V$ induces
the following structure on the subspace $V_{0}\oplus V_{1}$:
\begin{gather}
\b1\in V_0,\label{(2.2a)}
\\
T: \ \ V_{0}\rightarrow V_{1},\label{(2.2b)}
\\
{}_{(n)}:\ \ V_{i}\otimes V_{j}\rightarrow V_{(i+j-n-1)}, \quad n\geq 0,\quad
i,j=0,1,\label{(2.2c)}
\\
{}_{(-1)}: \quad (V_0\otimes V_{i})\oplus(V_{i}\otimes V_{0})\rightarrow
V_{i},\quad i=0,1.\label{(2.2d)}
\end{gather}
These data satisfy a list of identities obtained by inspecting those
listed in Def\/initions~2.2
and~2.6 
and choosing
the ones that make sense. The meaning of `make sense'  is clear in
the case of identities, such as the Jacobi, involving only
operations ${}_{(n)}$ with $n\geq 0$, because $V_0\oplus V_1$ is
closed under these operations. Expressions  $T^{j}(a_{(n)}b)$,
$a_{(i)}b_{(j)}c$, $(a_{(i)}b)_{(j)}c$, $a,b,c\in V_0,\;V_1$,
 are said to make sense if either they are compositions
of operations (\ref{(2.2a)})--(\ref{(2.2d)}) or because
$\Delta(a)+\Delta(b)-n-1<0$, $\Delta(b)+\Delta(c)-j-1<0$,
$\Delta(a)+\Delta(b)-i-1<0$ (resp.) in which case the expressions
are def\/ined to be 0, cf.~condition~(\ref{(2.1)}).
 Finally, we shall
say that an identity makes sense if all the expressions that it
involves make sense.

\paragraph{Def\/inition 3.2.} 
A {\it Courant
algebroid} is a vector space $V_{0}\oplus V_{1}$ carrying the data
(\ref{(2.2a)})--(\ref{(2.2d)}) so that all those axioms of Def\/inition~2.6
that make sense are valid.

\paragraph{Def\/inition 3.3.} 
 A {\it vertex  algebroid} is a vector space
$V_{0}\oplus V_{1}$ carrying the data (\ref{(2.2a)})--(\ref{(2.2d)})
so that all those axioms of Def\/inition~2.2
that make sense are valid.

\medskip

There are two obvious forgetful functors{\samepage
\begin{gather}
\Phi_{\rm alg}: \ \ \{\text{Vertex algebras}\}\rightarrow \{\text{Vertex
algebroids}\},\label{(2.3a)}
\\
\Phi^{\rm poiss}_{\rm alg}: \ \ \{\text{Vertex Poisson algebras}\}\rightarrow
\{\text{Courant algebroids}\},\label{(2.3b)}
\end{gather}
and both af\/ford the left adjoints $\CU_{\rm alg}$ and
$\CU^{\rm poiss}_{\rm alg}$, which are analogous to (\ref{(1.3b)}).}

When applied to the algebras of (\ref{(1.7)}) and of Examples~2.4, 2.5, they give us f\/irst examples of
Courant
\begin{equation}
\CV^{\rm poiss}(\Vir)=\BC,\qquad \CV^{\rm poiss}(\fg)=\BC\oplus\fg, \qquad \CV^{\rm poiss}(gl_N)=\BC\oplus gl_N  \label{(2.4)}
\end{equation}
and vertex algebroids  resp.
\begin{equation}
\CV(\Vir)_{c}=\BC,\qquad \CV(\fg)_{k}=\BC\oplus\fg,\qquad \CV(gl_N)_{k_1,k_2}=\BC\oplus
gl_N. \label{(2.5)}
\end{equation}
Here is an example of geometric nature.

\paragraph{Example 3.4.} 
 Let $A$ be a
commutative associative algebra with unit $\b1$, $\Omega(A)$ the
module of K\"ahler dif\/ferentials, $\Der(A)$ the algebra of
derivations. Note that the homological degree grading of $A$
determines that of $\Omega(A)$ and $\Der(A)$. Set $V_0=A$,
$V_1=\Der(A)\oplus\Omega(A)$, and let
$T\stackrel{\text{def}}{=}d:V_0=A\rightarrow\Omega_A\hookrightarrow
V_1$ be the canonical (de Rham) derivation. Then the space
$V_0\oplus V_1$ carries a unique Courant algebroid structure
determined by
\begin{gather}
 a_{(-1)}b=ab,\qquad a_{(-1)}(\tau+\omega)=a\tau+a\omega,\qquad
\tau_{(0)}a=\tau(a),\nonumber\\
\tau_{(0)}\omega=\text{Lie}_{\tau}\omega,\qquad \tau_{(1)}\omega=\iota_{\tau}\omega,\label{(2.6)}\\
\tau_{(0)}\xi=[\tau,\xi],\,\tau_{(1)}\xi=0,\nonumber
\end{gather}
for all $a,b\in A$, $\omega\in\Omega(A)$, $\tau,\xi\in \Der(A)$.
Denote thus def\/ined Courant algebroid by $\CV^{\rm poiss}(A)$.

\medskip

Since all the operations recorded in (\ref{(2.6)}) are of geometric
nature, there arises, for each scheme $X$, a sheaf of vertex
algebroids
\begin{equation}
\CV^{\rm poiss}_{X}\stackrel{\text{def}}{=}\CT_X\oplus\Omega_X.\label{(2.7)}
\end{equation}
Note that in keeping with our convention we assume that $X$ is {\it
graded}, that is to say, $\CO_X$ is a~sheaf of graded algebras;
consequently, $\CT_X$ and $\Omega_X$ are sheaves of graded
$\CO_X$-modules.

If $\CV=\CV_0\oplus\CV_1$ is a vertex algebroid, then part of its
structure coincides with that of $\CV^{\rm poiss}(A)$. For example, the
triple $(\CV_0,_{(-1)},\b1)$ is a commutative associative algebra
with unit, $(\CV_0)_{(-1)}(T(\CV_0))$ is a $\CV_0$-module (although
$\CV_1$ is not), the map $T:\CV_0\rightarrow
(\CV_0)_{(-1)}(T(\CV_0))$ is a~derivation. There arises a f\/iltration
\begin{equation}
\CV_0\oplus (\CV_0)_{(-1)}(T(\CV_0))\subset \CV,\label{(2.8)}
\end{equation}
and a moment's thought will show that, absolutely analogously to
Section~2.7, 
the corresponding ${\rm Gr}\,\CV$ carries a canonical
Courant algebroid structure.

\paragraph{Def\/inition 3.5.} \quad {} 

\begin{enumerate}\itemsep=0pt
\item[(a)] Given a Courant algebroid $\CV^{\rm poiss}$, call a vertex
algebroid $\CV$ a {\it quantization} of $\CV^{\rm poiss}$ if there is an
isomorphism of Courant algebroids
\[
{\rm Gr}\CV\iso\CV^{\rm poiss}.
\]
\item[(b)] If $\CV$ is a quantization of $\CV^{\rm poiss}(A)$, see
Example~3.4, 
then we shall denote $\CV$ by $\CV(A)$.
\end{enumerate}

Note that in the case of $\CV^{\rm poiss}(A)$, f\/iltration \eqref{(2.8)} implies
the following exact sequence of vector spaces
\begin{equation}
0\rightarrow\Omega(A)\rightarrow\CV_{1}\stackrel{\pi}{\rightarrow}
\Der(A)\rightarrow 0.\label{(2.9)}
\end{equation}

The problem of quantizing $\CV^{\rm poiss}(A)$ is not trivial and was
studied in \cite{GMSI}. Call $A$ {\it suitable for chiralization} if
$\Der(A)$ is a free $A$-module on generators $\tau_1,\dots ,\tau_N$ s.t.\
$[\tau_i,\tau_j]=0$. If $A$ is suitable for chiralization and a
basis $\tau_1,\dots ,\tau_N$ is f\/ixed, then $\CV^{\rm poiss}(A)$ can be
quantized by letting $\CV(A)=\CV^{\rm poiss}(A)$ as a vector space and
requiring all of the relations (\ref{(2.6)}) except the last two;
the latter are to hold only for the basis vector f\/ields:
\begin{equation}
(\tau_i)_{(0)}(\tau_j)=0,\qquad (\tau_i)_{(1)}(\tau_j)=0.\label{(2.10)}
\end{equation}
It is easy to show, using axioms (\ref{(1.1a)})--(\ref{(1.1e)}), at
least that these choices determine a vertex algebroid structure.

Let us  point out  the dif\/ferences between the Courant and vertex
algebroid structures thus obtained:  in the Courant case operation
$_{(1)}$ is an $A$-bilinear pairing and
\[
(f_{(-1)}\tau_i)_{(1)}(g_{(-1)}\tau_j)=0,
\]
while in the vertex case
\begin{equation}
(f_{(-1)}\tau_i)_{(1)}(g_{(-1)}\tau_j)=-f_{(-1)}(\tau_j(\tau_i(g))-g_{(-1)}(\tau_i(\tau_j(f))-
\tau_i(g)_{(-1)}\tau_j(f);\label{(2.11a)}
\end{equation}
in the Courant case multiplication $_{(-1)}$ is associative, e.g.,
\[
(fg)_{(-1)}\tau_{i}=f_{(-1)}(g_{(-1)}\tau_i),
\]
in the vertex case it is not as
\begin{equation}
(fg)_{(-1)}\tau_{i}=f_{(-1)}(g_{(-1)}\tau_i)-\tau_i(g)df
-\tau_i(f)dg;\label{(2.11b)}
\end{equation}
here $f,g\in A$ and axioms (\ref{(1.1c)})--(\ref{(1.1e)}) along with simplif\/ications due to grading have been
used.

Furthermore, if $A$ is suitable for chiralization, then
\begin{gather}
 \{\text{isomorphism classes of quantizations of }\CV^{\rm poiss}(A)\}\nonumber\\
 \qquad \qquad \text{is an} \quad
\big(\Omega^{3,{\rm cl}}(A)/d\Omega^{2}(A)\big)^{0}-\text{torsor},
\label{(2.12)}
\end{gather}
and
\begin{equation}
\text{Aut}(\CV(A))\iso\big(\Omega^{2,{\rm cl}}(A)\big)^{0},\label{(2.13)}
\end{equation}
where the automorphism of $\CV(A)$ attached to
$\omega\in\Omega^{2,{\rm cl}}(A)^0$ is def\/ined by the assignment
\begin{equation}
\CV(A)_{1}\rightarrow\CV(A)_{1},\qquad
\tau\mapsto\tau+\omega(\pi(\tau)),\label{(2.14)}
\end{equation}
$\pi$ being def\/ined in (\ref{(2.9)}). Note that we had to pay the
price for relentlessly working in the graded setting by extracting
the homological degree 0 subspace in (\ref{(2.12)}), (\ref{(2.13)}).

 Since any smooth algebraic variety $X$ can be covered by the
spectra of rings suitable for chiralization, (\ref{(2.12)}), (\ref{(2.13)}) create an avenue to def\/ine sheaves of vertex
algebroids over $X$, to be denoted $\CV_X$ or $\CV_A$ if
$X={\rm Spec}(A)$.

\paragraph{Def\/inition 3.6.} 
Let $X$ be a graded
scheme.

\begin{enumerate}\itemsep=0pt
\item[(a)] Call a sheaf of vertex algebroids over $X$ a {\it quantization}
of $\CV^{\rm poiss}_X$, see \eqref{(2.7)},  if for each af\/f\/ine $U\subset X$ its
space space of sections over $U$ is a quantization of
$\Gamma(U,\CV^{\rm poiss}_X)$.

\item[(b)] Denote by $\Vert_X$ (or $\Vert_A$ if $X={\rm Spec}(A)$) the category
of  quantizations of the sheaf $\CV_X^{\rm poiss}$.
\end{enumerate}

The characteristic property of $\CV_X\in\Vert_X$ is the existence of
the  sequence of sheaves
\begin{equation}
0\rightarrow\Omega_X\rightarrow\CV_X\stackrel{\pi}{\rightarrow}\CT_X\rightarrow
0; \label{(2.15)}
\end{equation}
this is  a sheaf analogue of (\ref{(2.9)}).

Here are  the (obvious graded versions of the) main results of
\cite{GMSI}:

\begin{itemize}\itemsep=0pt
\item there is a gerbe of vertex algebroids over a manifold $X$ such that the space of sections
over each $U={\rm Spec}(A)\subset X$, $A$ being suitable for
chiralization, is the category of quantizations of $\CV^{\rm poiss}(A)$;

\item in the case where $\CO_X=\CO_X^0$ this gerbe possesses a global section, i.e., a sheaf
$\CV_{X}\in\Vert_X$, if and only if the 2nd component of the Chern
character
\begin{equation}
{\rm ch}_{2}(\CT_{X})\in
H^{2}\big(X,\Omega^{2}_{X}\rightarrow\Omega^{3,{\rm cl}}_{X}\big)  \label{(2.16)}
\end{equation}
vanishes; if this class vanishes, then
\begin{equation}
\{\text{Isomorphism classes of sheaves }\CV_{X}\}  \ \ \text{is an} \ \
H^{1}\big(X,\Omega^{2}_{X}\rightarrow\Omega^{3,{\rm cl}}_{X}\big)-\text{torsor};\label{(2.17)}
\end{equation}
\item the forgetful functor (\ref{(2.3a)}) has a left adjoint functor
\begin{equation}
\CU_{\rm alg}: \ \ \{\text{Vertex algebroids}\}\rightarrow \{\text{Vertex
algebras}\};\label{(2.18)}
\end{equation}
if the obstruction (\ref{(2.16)}) vanishes, we call
\begin{equation}
D^{\rm ch}_{X}\stackrel{\text{def}}{=}\CU_{\rm alg}\CV_{X} \label{(2.19)}
\end{equation}
a sheaf of {\it chiral differential operators}, CDO for short.
\end{itemize}

 Note that proving (\ref{(2.17)})
amounts to covering $X$ by open sets that are suitable for
chiralization and re-gluing a given sheaf by composing the old
gluing functions with automorphisms (\ref{(2.13)}), (\ref{(2.14)}).

\paragraph{3.7.~Further examples and
constructions.} \quad {} 

\paragraph{3.7.1. Localization.} 
For any
quantization $\CV(A)$ and an ideal $\fa\subset A$ a natural
quantiza\-tion~$\CV(A_{\fa})$ is def\/ined \cite{GMSI}, the reason being
that all the $_{(n)}$-products on $\CV(A)$ are in fact certain
dif\/ferential operators. For example, at the quasiclassical level,
all the operations recorded in~(\ref{(2.6)}) are dif\/ferential
operators of order $\leq 1$. Furthermore, (\ref{(2.11a)}), (\ref{(2.11b)}) provide examples of genuine quantum operations being
order $\leq 2$ dif\/ferential operators. Therefore, given an $A$ and~$\CV(A)$, there arises a sheaf of vertex algebroids
\begin{equation}
\CV_{A}\in\Vert_A\quad \text{s.t.} \quad \Gamma({\rm Spec}(A),\CV_A)=\CV(A).\label{(2.20)}
\end{equation}
This construction underlies the above discussion of gerbes of vertex
algebroids.

A little more generally, if  $X$ and $Y$ are manifolds and
$p:X\rightarrow Y$ is a covering, then there is a functor
\[
p^*: \ \ \Vert_Y\rightarrow\Vert_X.
\]
The reason for this to be true is that the story of quantizing
$\CV^{\rm poiss}(A)$, ${\rm Spec}(A)\subset Y$, starts with a choice of an
Abelian basis $\{\tau_i\}\subset \Der(A)$, and any such choice is
canonically lifted to an Abelian basis of $\Gamma(V,\CT_X)$ for any
af\/f\/ine $V\subset p^{-1}({\rm Spec}(A))$.

In particular, if $X$ carries a free action of a f\/inite group $G$,
then there arises an equivalence of categories
\begin{equation}
\Vert_{X/G}\rightarrow\Vert_{X}^{G},\label{(2.21)}
\end{equation}
where $\Vert_{X}^{G}$ is a full subcategory of $G$-equivariant
vertex algebroids. The inverse functor is, essentially,  that of
$G$-invariants
\begin{equation}
\Vert_{X}^{G}\rightarrow\Vert_{X/G},\qquad \CV\mapsto
p_{*}(\CV^G),\label{(2.22)}
\end{equation}
where $p_*$ is the push-forward in the category of sheaves of vector
spaces.

Here is a version of localization called a {\it push-out} in
\cite{GMSI,GMSII}. If a Lie group acts on $A$ by derivations, then
\begin{equation}
A\otimes\CV(\fg)_k\stackrel{\rm def}{=}A\oplus
A\otimes\fg\oplus\Omega(A)\label{(2.23)}
\end{equation}
is a vertex algebroid.

\paragraph{3.7.2.~Af\/f\/ine space.} 
According to
(\ref{(2.16)}), (\ref{(2.17)}), there is a unique sheaf of vertex
algebroids over $\BC^N$, $\CV_{\BC^N}$. Its space of global sections
is $\CV(\BC[x_1,\dots ,x_N])$, where we take
$\{\partial_j=\partial/\partial x_j\}$ for an Abelian basis of
$\Der(\BC[x_1,\dots ,x_N])$, and if $U=\{f\neq 0\}$, then
$\Gamma(U,\CV_{\BC_N})$ is def\/ined via the localization of
Section~3.7.1. 

The corresponding CDO, see (\ref{(2.19)}), is
$D^{\rm ch}_{\BC^N}=\CU_{\rm alg}\CV_{\BC^N}$. Its space of global sections,
denoted by either $D^{\rm ch}(\BC^N)$ or $D^{\rm ch}(\BC[x_1,\dots ,x_N])$, is
likewise obtained via $\CU_{\rm alg}$:
$D^{\rm ch}(\BC^N)=\CU_{\rm alg}\CV(\BC[x_1,\dots ,x_N])$. More explicitly,
$D^{\rm ch}(\BC^N)$ can be def\/ined to be the vertex algebra gene\-ra\-ted
freely by the vector space $\BC^N\oplus(\BC^N)^*$ and relations
\[
(\partial_i)_{(0)}x_j=-(x_j)_{(0)}\partial_i=\delta_{ij}\b1,\qquad
a_{(n)}b=0\quad \text{for all} \quad a,b\in \BC^N\oplus(\BC^N)^*,  \quad n>0,
\]
where $\{\partial_i\}$, $\{x_j\}$ are dual bases of $\BC^N$ and
$(\BC^N)^*$ resp.

\paragraph{3.7.3.~Punctured plane.} 
Let
$X=\BC^{2}\setminus(0,0)$ and choose the trivial grading where
$\CO_X=\CO_X^0$. The obstruction (\ref{(2.16)}) vanishes, because
one sheaf, say, the restriction of $\CV_{\BC^{2}}$ to
$\BC^{2}\setminus(0,0)$, exists. Isomorphism classes of sheaves of
vertex algebroids over $X$ are easy to classify. Indeed, consider
the af\/f\/ine covering
\[
X=U_1\cup U_2,\qquad U_j=\{(y_1,y_2) \ \text{s.t.} \ y_j\neq 0\}
\]
and the Cech 1-cocycle
\[
\omega_{ab}: \ \ U_1\cap U_2\mapsto\frac{dy_1\wedge
dy_2}{y_{1}^{a}y_{2}^{b}},\quad a,b\geq 1.
\]
It is known (and easy to check) that
\[
H^{1}\big(X,\Omega^{2}_{X}\rightarrow\Omega^{3,{\rm cl}}_{X}\big)=\oplus_{a,b}\BC\omega_{ab}.
\]
Hence the isomorphism classes of sheaves of vertex algebroids over
$X$ are in 1-1 correspondence with linear combinations of
$\omega_{ab}$. Here is an explicit construction of the sheaf
attached to $k\omega_{ab}$: let  $\CV_{X}$ be the restriction of the
standard $\CV_{\BC^{2}}$ to $X$, $\CV_{U_{j}}$ its pull-back  to
$U_j$, $j=1,2$; now glue~$\CV_{U_{1}}$ and~$\CV_{U_{2}}$ over the
intersection $U_1\cap U_2$, cf.~(\ref{(2.13)}), (\ref{(2.14)}), as
follows:
\begin{equation}
\partial_{1}\rightarrow\partial_{1}+\frac{kT(y_2)}{y_{1}^{a}y_{2}^{b}},\qquad
\partial_{2}\rightarrow\partial_{2}-\frac{kT(y_1)}{y_{1}^{a}y_{2}^{b}}.
\label{(2.24)}
\end{equation}
It is immediate to generalize this to the case of an arbitrary~$\omega$ in the linear span of~$\{\omega_{ab}\}$. Denote the sheaf
sheaf thus def\/ined by $\CV_{X,\omega}$.

 This simple
example will be essential for our purposes.

\medskip

\noindent
{\bf 3.7.4. Symmetries.} 
Given a Lie
algebra morphism
\[
\rho: \ \ \fg\rightarrow \Der(A),
\]
the composition
\[
\fg\rightarrow
\Der(A)\stackrel{(\text{id},0)}{\hookrightarrow}\CV^{\rm poiss}(A)=\Der(A)\oplus\Omega(A)
\]
def\/ines a Courant algebroid morphism
\[
\rho^{\rm poiss}: \ \ \CV^{\rm poiss}(\fg)\rightarrow\CV^{\rm poiss}(A).
\]
Furthermore,  any Poisson algebroid morphism $\rho^{\rm poiss}$, upon
quotienting out by $\Omega(A)$,  def\/ines a~Lie algebra morphism
$\rho$.

\medskip

A {\it quantization} of  a Courant algebroid morphism
$\rho^{\rm poiss}:\CV^{\rm poiss}(\fg)\rightarrow\Gamma(X,\CV^{\rm poiss}_{X})$
is def\/ined to be a vertex algebroid morphism
$\hat{\rho}:\CV(\fg)_{k}\rightarrow\CV(A)$ such that the diagram
\begin{equation}
\xymatrix{ 0\ar[r]&\Omega(A)\ar[r]&\CV(A)\ar[r]&
\Der(A)\ar[r]&0\\
 & & \CV(\fg)_{k}\ar[r]^{\rho}\ar[u]^{\hat{\rho}}&\Der(A)\ar@{=}[u]\ar[r] & 0
 }
 \label{(2.25)}
\end{equation}
commutes. Here the arrow $\CV(\fg)_{k}\stackrel{\rho}{\rightarrow}
\Der(A)$ means the composition
$\CV(\fg)_{k}=\BC\oplus\fg\stackrel{0\oplus\rho}{\rightarrow}
\Der(A)$, see (\ref{(2.5)}).

Similarly, if $\fg$ operates on  $X$, that is, there is a Lie
algebra morphism
\[
\rho: \ \ \fg\rightarrow \Gamma(X,\CT_{X})
\]
then  a {\it quantization} of the corresponding Courant algebroid
morphism
$\rho^{\rm poiss}:\CV^{\rm poiss}(\fg)\rightarrow\Gamma(X,\CV^{\rm poiss}_{X})$
is def\/ined to be a vertex algebroid morphism
$\hat{\rho}:\CV(\fg)_{k}\rightarrow\Gamma(X,\CV_{X})$ such that the
diagram
\begin{equation}
\xymatrix{
0\ar[r]&\Gamma(X,\Omega_{X})\ar[r]&\Gamma(X,\CV_{X})\ar[r]&
\Gamma(X,\CT_{X})\ar[r]&0\\
 & & \CV(\fg)_{k}\ar[r]^{\rho}\ar[u]^{\hat{\rho}}&\Gamma(X,\CT_{X})\ar@{=}[u]\ar[r]&0
 }
 \label{(2.26)}
\end{equation}
commutes.

To see an example of importance for what follows, let us consider
the tautological action of~$gl_2$ on $\BC^2$. If we let
$X=\BC^2\setminus(0,0)$, then there arises
\begin{equation}
\rho^{\rm poiss}: \ \ \CV^{\rm poiss}(gl_2)\rightarrow
\Gamma\big(\BC^2\setminus(0,0),\CV^{\rm poiss}_{\BC^2\setminus(0,0)}\big),\label{(2.27)}
\end{equation}
and we ask if this map can be quantized to a map
$\CV(gl_2)_{k_1,k_2}\rightarrow
\Gamma(\BC^2\setminus(0,0),\CV_{\BC^2\setminus(0,0),\omega})$, where
the vertex algebroid $\CV_{\BC^2\setminus(0,0),\omega}$ was def\/ined
in Section~3.7.3. 

\paragraph{Lemma 3.7.5.} 
{\it Quantization of \eqref{(2.27)},
\[
\hat{\rho}: \ \ \CV(gl_2)_{k_1,k_2}\rightarrow
\Gamma\big(\BC^2\setminus(0,0),\CV_{\BC^2\setminus(0,0),\omega}\big)
\]
exists if and only if $\omega=kdy_1\wedge dy_2/y_1 y_2$, $k_1=-k-1$,
$k_2=k-1$ for some $k\in\BC$.}


\begin{proof} We shall use  the notation of Section~3.7.3. 
In terms of the coordinates $y_1$, $y_2$ the morphism $\rho$ is this
\begin{equation}
\rho(E_{ij})\mapsto y_i\partial_j.\label{(2.28)}
\end{equation}
If we consider $y_i\partial_j$ as an element of $\Gamma(U_1,\CV_{X,\omega})$, then over $U_2$
it becomes, according to (\ref{(2.24)}), $(y_i\partial_j\pm kT(y_{j\pm 1})y_i)/y_1^a y_2^b$ and
hence may develop a pole. To compensate for it, we can choose a dif\/ferent lift of $E_{ij}$ to
$\Gamma(U_1,\CV_{X,\omega})$ by replacing $y_i\partial_j$ with $y_i\partial_j+\alpha$, where
$\alpha\in \Gamma(U_1,\Omega_X)$. Over $U_2$ this element becomes $(y_i\partial_j\pm T(y_{j\pm
1})y_i)/y_1^a y_2^b+\alpha$. Since $\alpha$ may have at most a pole along $\{y_1=0\}$, for this
element to extend to a section over $U_2$, one of the following two things must happen: either
$T(y_{j\pm 1})y_i/y_1^a y_2^b$ has no pole along $\{y_1=0\}$, in which case no $\alpha$ is
needed, or $T(y_{j\pm 1})y_i/y_1^a y_2^b$ has no pole along $\{y_2=0\}$, in which case a
desired $\alpha$ can be found. For a~favorable event to occur for $i=1$ and $i=2$, both $a$ and
$b$ must be at most~1. But by def\/inition, see Section~3.7.3, 
$a$ and $b$ are at least~1; therefore a linear map $gl_2\rightarrow \Gamma(X,\CV_{X,\omega})$ may exist only if
$\omega=k dy_1\wedge dy_2/y_1y_2$. On the other hand, if $\omega=k dy_1\wedge dy_2/y_1y_2$,
then the map $\hat{\rho}$ def\/ined so that
\begin{gather}
 \hat{\rho}(E_{12})= y_1\partial_{y_{2}},\label{(2.29a)}
\\
 \hat{\rho}(E_{21})= y_2\partial_{y_{1}}-\frac{ky_{2}'}{y_1},\label{(2.29b)}
\\
\hat{\rho}(E_{11})= y_{1}\partial_{y_{1}},\label{(2.29c)}
\\
 \hat{\rho}(E_{22})= y_2\partial_{y_{2}}+\frac{ky_{1}'}{y_1}.\label{(2.29d)}
 \end{gather}
delivers the desired vertex algebroid morphism
\begin{equation}
\hat{\rho}:\ \ \CV(gl_2)_{-k-1,k-1}\rightarrow\Gamma\big(\BC^2\setminus(0,0),\CV_{\BC^2\setminus(0,0),\omega}\big), \qquad \omega=k
dy_1\wedge dy_2/y_1y_2.\label{(2.30)}
\end{equation}
It is easy to see the vertex algebroid morphism condition determines
the map uniquely.
\end{proof}

Note that the top row of (\ref{(2.26)}), unlike that of
(\ref{(2.25)}), does not have to be exact in general. In the case at
hand, however, it is precisely when $\omega=k\omega_{11}$:

\paragraph{Corollary 3.7.6.} 
{\it The
sequence
\begin{gather*}
0\rightarrow\Gamma\big(\BC^2\setminus(0,0),\Omega_{\BC^2\setminus(0,0)}\big)\rightarrow
\Gamma\big(\BC^2\setminus(0,0),\CV_{\BC^2\setminus(0,0),\omega}\big)_{1}\rightarrow
\Gamma\big(\BC^2\setminus(0,0),\CT_{\BC^2\setminus(0,0)}\big)\rightarrow 0
\end{gather*}
is exact if and only if $\omega=k\omega_{11}$ for some $k\in\BC$.}

\begin{proof}
 Notice that $\CT_{\BC^2\setminus(0,0)}$ is generated by
$\rho(gl_2)$ over functions. The ``if'' part is then seen to be an
immediate consequence of Lemma~3.7.5. 
The ``only if''
part was actually proved at the beginning of the proof of the lemma
cited.
\end{proof}

\paragraph{3.7.7.~Conformal structure.}  
If
$x_1,\dots ,x_N$ are coordinates on $\BC^N$,
$\partial_j=\partial/\partial_j$, then there is a vertex (Poisson)
algebra morphism
\begin{gather}
\CU^{\rm poiss}\CL(\Vir)\rightarrow\Gamma\big(\BC^{N},\CD^{\rm poiss}_{\BC^N}\big),\qquad
 \CU\CL(\Vir)_{2N}\rightarrow\Gamma\big(\BC^{N},\CD_{\BC^N}\big),\nonumber\\
 L\mapsto\sum_{j=1}^{N}T(x_j)_{(-1)}\partial_{j},\label{(2.31)}
\end{gather}
where the latter is a quantization of the former.  A little more
generally, if $A$ is an algebra suitable for chiralization with
$\tau_1,\dots ,\tau_N$ an Abelian basis of $\Der(A)$, then one can f\/ind
a coordinate system, i.e., $\{x_1,\dots ,x_N\}\subset A$ s.t.
$\tau_i(x_j)=\delta_{ij}$, and thus obtain
\begin{gather*}
\CU^{\rm poiss}\CL(\Vir)\rightarrow\CU^{\rm poiss}\CV^{\rm poiss}(A),\qquad
 \CU\CL(\Vir)_{2N}\rightarrow\CU\CV(A),\qquad L\mapsto\sum_{j=1}^{N}T(x_j)_{(-1)}\tau_{j}.
\end{gather*}
In this case, $N$ is the Krull dimension of $A$.

Another example is provided by the twisted sheaves
$\CV_{\BC^2\setminus 0,\omega}$ of Section~3.7.3. 
Somewhat
unexpectedly, the same def\/inition (\ref{(2.31)}), which in the
present case becomes $L\mapsto T(y_1)\partial_1+T(y_2)\partial_2$,
applied locally over both both  charts $U_1$ and $U_2$ survives the
twisted gluing transformation (\ref{(2.24)}) for any $\omega$ and
def\/ines a global morphism
\begin{equation}
\CU\CL(\Vir)_2\rightarrow\Gamma\big(\BC^2\setminus 0,\CV_{\BC^2\setminus
0,\omega}\big).\label{(2.32)}
\end{equation}

 \section{A graded Courant algebroid\\
  attached to a commutative associative algebra}
\label{sect.3}

 \paragraph{4.1.~Modules of dif\/ferentials.} 
Even though the assumption that all the vector spaces in question are $\BZ$-graded has been
 kept since the very beginning of Section~\ref{sect.1}, it has  been barely used. From now on it will be essential
 and referred to explicitly as in the following def\/inition.

\paragraph{Def\/inition 4.1.1.} 
A {\it differential graded algebra} (DGA) $R$ is a pair
$(R_{\#},D)$,
 where $R_{\#}=\oplus_{n=0}^{\infty}R^{n}$ is
 a graded supercommutative associative algebra with $R^{\rm even}=\oplus_{n=0}^{\infty}V^{2n}$,
 $R^{\rm odd}=\oplus_{n=0}^{\infty}R^{2n+1}$, and $D$ is a square 0 degree $(-1)$ (hence odd)
 derivation. Call a DGA $(R_{\#},D)$ {\it quasi-free} if there is a graded vector superspace
  $V=\oplus_{n=0}^{\infty}V^{n}$
  with $V^{\rm even}=\oplus_{n=0}^{\infty}V^{2n}$, $V^{\rm odd}=\oplus_{n=0}^{\infty}V^{2n+1}$
 such that $R_{\#}$ is the symmetric algebra $S^{\bullet}V$.

 \medskip

If $R$ is a DGA, then $H^{\bullet}_{D}(R)\stackrel{\text{def}}{=}\Ker(D)/{\rm Im}(D)$ is a a graded
supercommutative associative algebra.

For any commutative associative algebra $A$ there is a quasi-free DGA $R$ and a
quasi-iso\-mor\-phism
\begin{equation}
R\rightarrow A,\label{(3.1.1)}
\end{equation}
that is to say, a DGA morphism ($A$ being placed in homological degree 0 and equipped with a~zero dif\/ferential) that delivers a graded algebra isomorphism $H^{\bullet}_{D}(R)\iso A$.

If $A$ is f\/initely generated, then a DGA resolution $R$ can be chosen so that each $V^j$ from
Def\/inition~4.1.1 
is f\/inite dimensional. {\it These two finiteness assumptions will
be made throughout.}

 A DGA resolution of $A$ is not unique, but for any two such resolutions
\[
R_1\rightarrow A\leftarrow R_2
\]
there is a homotopy equivalence \cite{Behr}
\begin{equation}
f: \ \ R_1\rightarrow R_2.\label{(3.1.2)}
\end{equation}

If $R$ is a quasi-free DGA, denote by $\Omega(R)$  the module of K\"ahler dif\/ferentials of $R$.
It is cano\-ni\-cally a dif\/ferential-graded (DG) free $R$-module with derivation $d: R\rightarrow
\Omega(R)$ and dif\/ferential ${\rm Lie}_{D}$, which we choose to denote by $D$, too. The
correspondence $R\mapsto\Omega(R)$ is functorial in that naturally associated to an algebra
morphism $f:R_1\rightarrow R_2$ there is a map of DG $R_1$-modules:
\begin{equation}
\Omega(f): \ \ \Omega(R_1)\rightarrow \Omega(R_2)\label{(3.1.3)}.
\end{equation}
Furthermore, we have
\begin{equation}
\Omega(R)=\bigoplus_{n=0}^{+\infty} \Omega(R)^{n},\qquad  d: \ R^{n}\rightarrow \Omega(R)^{n},\quad
D:\ \Omega(R)^{n}\rightarrow \Omega(R)^{n-1},\quad  [d,D]=0.\label{(3.1.4)}
\end{equation}
It follows that the homology $H_{D}(\Omega(R))$ is naturally a graded
$H^{\bullet}_{D}(R)$-module.

For any 2 quasi-free DGA resolutions $R_1\rightarrow A\leftarrow R_2$, we can f\/ind a homotopy
equivalence $f: R_1\rightarrow R_2$, see (\ref{(3.1.2)}), hence a quasi-isomorphism
\begin{equation}
\Omega(f): \ \ \Omega(R_1)\rightarrow \Omega(R_2) \label{(3.1.5)}
\end{equation}
and an isomorphism
\begin{equation}
H(\Omega(f)): \ \ H_{D_1}^{\bullet}(\Omega(R_1))\rightarrow H_{D_2}^{\bullet}(\Omega(R_2)).
\label{(3.1.6)}
\end{equation}

\paragraph{Def\/inition 4.1.2.} 
\begin{equation}
\Omega(A)^{\bullet}=H_{D}^{\bullet}(\Omega(R)),\label{(3.1.7)}
\end{equation}
where $R$ is a quasi-free DGA resolution of $A$.

\medskip

The assignment $A\mapsto\Omega(A)^{\bullet}$ def\/ines a functor from the category of algebras to
the category of graded vector spaces.

Note that
\begin{equation}
\Omega(A)^{\bullet}=\bigoplus_{n=0}^{+\infty}\Omega(A)^{n},\label{(3.1.8)}
\end{equation}
is a graded $A$-module, and $\Omega(A)^{0}$ is the module of K\"ahler dif\/ferentials of $A$,
$\Omega(A)$.

\paragraph{4.2. Modules of derivations.} 
If $R$ is a quasi-free DGA, we denote by $\Der(R)$ the Lie algebra of derivations of $R$. Like
$\Omega(R)$, it is a DG $R$-module, but unlike $\Omega(R)$ it is graded in both directions:
\begin{equation}
\Der(R)=\bigoplus_{n\in\BZ}\Der(R)^{n}\label{(3.2.1)}
\end{equation}
and, which is  more serious, not free; in fact, each component $\Der(R)^{n}$ is a direct product
\begin{equation}
\Der(R)^{n}=\prod_{j=0}^{+\infty}(V^{j})^{*}\otimes R^{n+j},\label{(3.2.2)}
\end{equation}
where $V^j$ is one of the ingredients of Def\/inition~4.1.1 
assumed to be f\/inite
dimensional.

The derivation $[D,\cdot]:\Der(R)\rightarrow \Der(R)$ is a dif\/ferential because $D\in \Der(R)^{-1}$ is
odd. Hence a Lie algebra $H_{[D,\cdot ]}^{\bullet}(\Der(R))$ arises.

The assignment $R\mapsto \Der(R)$ is not quite functorial, because even if $f: R_1\rightarrow
R_2$ is a quasi-isomorphism, a Lie algebra morphism $\Der(f):\Der(R_1)\rightarrow \Der(R_2)$ does
not quite exist. It does exist though at the level of the corresponding homotopy categories.
This remark and what follows belongs to Hinich \cite[Section~8]{Hin}.

Decompose $f: R_1\rightarrow R_2$ as follows
\begin{equation}
f: \ \ R_1\stackrel{i}{\hookrightarrow}S\stackrel{p}{\rightarrow}R_2,\label{(3.2.3)}
\end{equation}
where $i$ is a standard acyclic cof\/ibration, and $p$ is an acyclic f\/ibration. (Recall that $i$
being a~standard cof\/ibration means  $S$ being obtained by adjoining variables to $R_1$, and
being a f\/ibration means being an epimorphism.)

In the case of $i$, there arises a diagram of quasiisomorphisms
\begin{equation}
\Der(R_1)\stackrel{\pi_{i}}{\leftarrow} \Der(i)\stackrel{{\rm in}_{i}}{\rightarrow}
\Der(S),\label{(3.2.4)}
\end{equation}
where $\Der(i)=\{\tau\in \Der(S) \ \text{s.t.} \ \tau(R_1)\subset R_1\}$, $\pi_{i}$ is the obvious
projection, and ${\rm in}_{i}$ is the obvious embedding.

Analogously, in the case of $p$, there is a diagram of quasiisomorphisms
\begin{equation}
\Der(S)\stackrel{{\rm in}_{p}}{\leftarrow} \Der(p)\stackrel{\pi_{p}}{\rightarrow}
\Der(S),\label{(3.2.5)}
\end{equation}
where $\Der(p)=\{\tau\in \Der(R_1)\ \text{s.t.} \ \tau(\Ker(p))\subset \Ker(p)\}$, ${\rm in}_{p}$ is the
obvious embedding, and $\pi_{p}$ is the obvious projection.

Hinich def\/ines
\begin{equation}
\Der(i)={\rm in}_{i}\circ \pi_{i}^{-1},\qquad \Der(p)=\pi_{p}\circ {\rm in}_{p}^{-1}, \qquad \Der(f)=\Der(p)\circ \Der(i).
\label{(3.2.6)}
\end{equation}
This map makes sense in the homotopy category and delivers a homotopy category
quasi-isomorphism
\begin{equation}
\Der(f): \ \ \Der(R_1)\rightarrow \Der(R_2).\label{(3.2.7)}
\end{equation}
Hence an isomorphism
\begin{equation}
H(\Der(f)): \ \ H_{[D_1,\cdot ]}^{\bullet}(R_1)\rightarrow H_{[D_2,\cdot ]}^{\bullet}(R_2).\label{(3.2.8)}
\end{equation}

\paragraph{Theorem 4.2.1 (\cite{Hin}).} 
{\it If $f_j:R_1\rightarrow R_2$,
$j=1,2$, are homotopic to each other, then $\Der(f_j)$, $j=1,2$, are also.}

\paragraph{Corollary 4.2.2 (\cite{Hin}).} {\it \quad {} 
\begin{enumerate}\itemsep=0pt
 \item[$(i)$] The assignment $R\mapsto \Der(R)$ defines a functor from the homotopy category of DG
commutative associative  algebras with quasi-isomorphisms to the homotopy category of DG Lie
algebras.

\item[$(ii)$] The assignment $R\mapsto H^{\bullet}_{[D,\cdot ]}(\Der(R))$ defines a functor from the homotopy
category of DG commutative associative algebras with quasi-isomorphisms to the category of
graded Lie algebras.
\end{enumerate}}

\paragraph{Def\/inition 4.2.3.} 
\[
 \Der(A)^{\bullet}=H^{\bullet}_{[D,\cdot ]}(\Der(R)),
\]
where $R$ is a quasi-free DGA resolution of $A$.

\paragraph{4.3.~Synthesis: Courant algebroids.} 

The notion of a Courant algebroid allows us to bring Sections~4.1 
and~4.2 
under the same roof. Let $\CV^{\rm poiss}$ be a Courant algebroid. It follows from the Jacobi
identity (\ref{(1.1d)}) that, for any $\xi\in\CV^{\rm poiss}(R)$, $\xi_{(0)}\in
{\rm End}(\CV^{\rm poiss}(R))$ is a derivation of all products. Identity~(\ref{(1.1b)}) implies that
$\xi_{(0)}$ commutes with $T$. If, in addition, $\xi$ is odd and $\xi_{(0)}\xi=0$, then
$(\xi_{(0)})^{2}=0$ as another application of (\ref{(1.1d)}) shows. Therefore, a pair
$(\CV^{\rm poiss},\xi)$ is a dif\/ferential Courant algebroid, and the homology Courant algebroid,
$H^{\bullet}_{\xi_{(0)}}(\CV^{\rm poiss})$, arises.

 Let us now specialize this well-known construction to the
  Courant algebroid $\CV^{\rm poiss}(R)=\Der(R)\oplus\Omega(R)$, see
  Example~3.4, 
  in the case of a
 quasi-free DGA $R=(R_{\#},D)$.
By def\/inition $D$ is odd and, according to (\ref{(2.6)}), $D_{(0)}D=[D,D]=0$. Hence the pair
$(\CV^{\rm poiss}(R), D_{(0)})$ is a DG Courant algebroid and the graded Courant algebroid
$H_{D_{(0)}}^{\bullet}(\CV^{\rm poiss}(R))$ arises. Again by virtue of (\ref{(2.6)}), the
dif\/ferential $D_{(0)}$ preserves $\Der(R)\subset\CV^{\rm poiss}(R)$, where it coincides with
$[D,\cdot ]$, and~$\Omega(R)$, where it  coincides with the standard action of $D$ by the Lie
derivative, see also (\ref{(3.1.4)}). We obtain a canonical vector space isomorphism
\begin{equation}
H_{D_{(0)}}^{\bullet}(\CV^{\rm poiss}(R))=H_{[D,\cdot ]}^{\bullet}(\Der(R))\oplus
H_{D}^{\bullet}(\Omega(R)).\label{(3.3.1)}
\end{equation}
If $f:R_1\rightarrow R_2$ is a homotopy equivalence, then
\begin{equation}
H(\Der(f))\oplus H(\Omega(f)): \ \ H_{D_{(0)}}^{\bullet}(\CV^{\rm poiss}(R_1))\rightarrow
H_{D_{(0)}}^{\bullet}(\CV^{\rm poiss}(R_2)),\label{(3.3.2)}
\end{equation}
is a vector space isomorphism by virtue of (\ref{(3.1.6)}) and (\ref{(3.2.8)}). In fact,
(\ref{(3.3.2)}) is a graded Courant algebroid isomorphism. This follows from the fact that the
Courant algebroid structure on $\CV^{\rm poiss}(R)$ consists of classical dif\/ferential geometry
operations, such as the tautological action of~$\Der(R)$ on $R$ and the action of~$\Der(R)$ on
$\Omega(R)$ by means of the Lie derivative. An
 inspection of maps (\ref{(3.2.3)})--(\ref{(3.2.8)}) shows that Hinich's
construction respects all these operations.

\paragraph{Corollary 4.3.1.} 
{\it The assignment $R\mapsto
H^{\bullet}_{D_{(0)}}(\CV^{\rm poiss}(R))$ defines a functor from the homotopy cate\-gory of DG
commutative associative  algebras with quasi-isomorphisms to the category of graded Courant
algebroids.}

\paragraph{Def\/inition 4.3.2.} 
\[
 \CV^{\rm poiss}(A)^{\bullet}=H^{\bullet}_{D_{(0)}}(\CV^{\rm poiss}(R)),
\]
where $R$ is a quasi-free DGA resolution of $A$.

\paragraph{4.4. Conformal structure.} 
The construction of Section~3.7.7 
in the present setting means the following. If
$R_{\#}=S^{\bullet}V$, pick a homogeneous basis $\{x_i\}\subset V$ and a dual `basis'
$\{\partial_i\}\subset V^*$, where $\partial_i(x_j)=\delta_{ij}$. As in Section~3.7.7, 
we obtain a morphism
\begin{equation}
\CU^{\rm poiss}\CL(\Vir)\rightarrow \CU\CV^{\rm poiss}(R),\qquad L\mapsto\sum_{j}(T(x_j))_{(-1)}\partial_j.
\label{(3.4.1)}
\end{equation}

\paragraph{Lemma 4.4.1.} 
\begin{equation}
\xi_{(0)}\sum_{j}T(x_j)_{(-1)}\partial_j=0 \qquad \text{\it{for any}} \quad \xi\in \Der(R).\label{(3.4.2)}
\end{equation}

\paragraph{Corollary 4.4.2.} 
{\it Assignment \eqref{(3.4.1)} determines a Courant
algebroid morphism}
\begin{equation}
\CU^{\rm poiss}\CL(\Vir)\rightarrow \CU\CV^{\rm poiss}(A)^{0}. \label{(3.4.3)}
\end{equation}

\begin{proof}[Proof of Lemma.] Let $\xi=\sum_i (f_{i})_{(-1)}\partial_i$, $f_i\in R_{\#}$. We have
\begin{gather*}
 \xi_{(0)}\sum_{j}T(x_j)_{(-1)}\partial_j=\sum_{j}(\xi_{(0)}T(x_j))_{(-1)}\partial_j
 +\sum_{j}(-1)^{\xi\cdot x_j} T(x_j)_{(-1)}(\xi_{(0)}\partial_j) \\
\hphantom{\xi_{(0)}\sum_{j}T(x_j)_{(-1)}\partial_j}{} = \sum_{j}T(\xi_{(0)}x_j)_{(-1)}\partial_j
 -\!\sum_{j}(-1)^{\xi\cdot x_j+\xi\cdot x_j} T(x_j)_{(-1)}\!\left(\!(\partial_j)_{(0)}\!\sum_i (f_{i})_{(-1)}\partial_i\!\right) \!\! \\
\hphantom{\xi_{(0)}\sum_{j}T(x_j)_{(-1)}\partial_j}{}= \sum_{j}T(f_j)_{(-1)}\partial_j
 -\sum_{i}\sum_{j}T(x_j)_{(-1)}\left(\frac{\partial f_{i}}{\partial x_j}\right)_{(-1)}\partial_i
\\
\hphantom{\xi_{(0)}\sum_{j}T(x_j)_{(-1)}\partial_j}{} =\sum_{j}T(f_j)_{(-1)}\partial_j
 -\sum_{i}T(f_i)_{(-1)}\partial_i=0.\tag*{\qed}
 \end{gather*}\renewcommand{\qed}{}
\end{proof}

\paragraph{4.5. Grading of $\boldsymbol{\Der(A)^{\bullet}}$ and identif\/ication of
$\boldsymbol{\Der(A)^{0}}$.} 
Recall that $R$, hence $\Omega(R)$ and~$\Omega(A)^{\bullet}$, are all graded by~$\BZ_{+}$.
Contrary to this, although the Lie algebra $\Der(R)$ is graded in both directions, $\Der(A)$ is
$\BZ_{-}$-graded.

\paragraph{Lemma 4.5.1.} {\it 

\begin{enumerate}\itemsep=0pt
\item[$(a)$] $\Der(A)^n=0$ if $n\geq 0$;

\item[$(b)$] $\Der(A)^0$ is the Lie algebra of derivations of $A$.
\end{enumerate}}

\begin{proof} Consider a quasi-free DGA resolution $R\rightarrow A$. The complex $(\Der(R),
[D,\cdot ])$ is f\/iltered as follows, cf.~(\ref{(3.2.2)}),
\begin{equation}
F^{p}\, \Der(R)^{n}=\prod_{j=p}^{+\infty}(V^{j})^{*}\otimes R^{n+j}.\label{(3.5.1)}
\end{equation}
A spectral sequence $(E_{pq}^r, d_r)\Rightarrow \Der(A)^{p+q}$ arises so that
\[
(E^0_{pq},d_0)=\big((V^{-p})^{*}\otimes R^{q}, 1\otimes D\big).
\]
Since $R=(R_{\#},D)$ is quasi-isomorphic to $A$ placed in degree 0, we have
\[
E^1_{pq}=\left\{
\begin{array}{ll}
(V^{-p})^{*}\otimes A& \text{if} \ \ q=0,\\
0&\text{otherwise}.
\end{array}
\right.
\]
It follows at once that the spectral sequence collapses and $\Der(A)^{-n}$ is the $n$-th
cohomology of the complex
\begin{equation}
0\rightarrow (V^0)^{*}\otimes A\rightarrow (V^1)^{*}\otimes A\rightarrow\cdots\rightarrow
(V^n)^{*}\otimes A\rightarrow\cdots.\label{(3.5.2)}
\end{equation}
Item $(a)$ of the lemma is thus proven.

In order  to prove item $(b)$, we have to write down a formula for the dif\/ferential of complex~(\ref{(3.5.2)}). The resolution $R\rightarrow A$ gives an exact sequence
\[
0\rightarrow J\hookrightarrow R^{0}\stackrel{\pi}{\rightarrow} A\rightarrow 0.
\]
We shall regard an element $\tau\in (V^{j})^*$ as a derivation of $R=S^{\bullet}V$. The
dif\/ferential, $D$, of $R$  can be written thus: $D=\sum f_j\partial_{j}+\xi$, where
$\partial_{j}\in (V^{1})^*$, $\{f_j\}$ generate $J$, and $\xi\in F^{2} \Der(R)^{-1}$.
 It easily follows from the construction of the spectral
sequence that if $\tau\in (V^{0})^*$, then
\[
d^1(\tau\otimes a)=- a\pi(\tau(f_j))\partial_j.
\]
It follows at once that $\Ker\{d^1:(V^0)^{*}\otimes A\rightarrow (V^1)^{*}\otimes A\}$ is
precisely the algebra of derivations of $R^0$ that preserve the ideal $J$ modulo those
derivations whose image is $J$, and this is $\Der(A)^0$ by def\/inition.
\end{proof}

\section{Quantization in the case of a Veronese ring}\label{sect.4}

\paragraph{5.1.~Set-up.} 
Consider the Veronese ring
\begin{equation}
A_N=\BC[x_0,\dots ,x_N]/Q,\qquad Q=(x_i x_j-x_{i+1}x_{j-1}).\label{(4.1.0)}
\end{equation}
It is known that ${\rm Spec}(A_{N})$ is the cone over the highest weight vector orbit in the
projectivization of the $(N+1)$-dimensional representation of $sl_{2}$. Hence the canonical Lie
algebra morphism
\begin{equation}
sl_{2}\rightarrow \Der(A_N).\label{(4.1.1)}
\end{equation}
An explicit formula for this morphism will appear in Section~5.2.3 
below.

Being a cone, ${\rm Spec}(A_N)$ carries the Euler vector f\/ield $\sum_{j}x_{j}\partial_{j}$. This
allows us to extend (\ref{(4.1.1)}) to an action of $gl_2$:
\begin{equation}
gl_{2}\rightarrow \Der(A), \qquad \text{where} \qquad E_{11}+E_{22}\mapsto
N\sum_{j}x_{j}\partial_{j}.\label{(4.1.2)}
\end{equation}

\noindent
{\bf Remark.} The normalizing factor of $N$ is not particularly important but can be justif\/ied
by the geometry of the base af\/f\/ine space $SL_{2}/\CN$.

\medskip

As in Section~3.7.4, 
this gives  a~Courant algebroid morphism
\begin{equation}
\CV^{\rm poiss}(gl_2)\rightarrow\CV^{\rm poiss}(A_N)^{0}\subset
\CV^{\rm poiss}(A_N)^{\bullet}.\label{(4.1.3)}
\end{equation}

The following theorem, the main result of this paper, uses the concept of  quantization of a~Courant algebroid, see Def\/inition~3.5, 
and the notion of quantization of a~Courant
algebroid map, see Section~3.7.4, 
(\ref{(2.25)}), (\ref{(2.26)}).

\paragraph{Theorem 5.1.1.} {\it 
\begin{enumerate}\itemsep=0pt
\item[$(a)$] The Courant algebroid $\CV^{\rm poiss}(A_N)^{\bullet}$ admits a unique quantization
$\CV(A_N)^{\bullet}$.

\item[$(b)$] Maps \eqref{(4.1.2)} and \eqref{(3.4.3)} can be quantized to the maps
\begin{gather}
\CV(gl_2)_{-N-2,N}\rightarrow\CV(A_N)^{0}, \label{(4.1.4)}
\\
\CU\CL(\Vir)_{2}\rightarrow\CU_{\rm alg}\CV(A_N)^{0}, \label{(4.1.5)}
\end{gather}
where the functors $\CU$ and $\CU_{\rm alg}$ are those defined in \eqref{(1.3b)}  and \eqref{(2.18)}
resp.
\end{enumerate}}

The proof of the theorem is constructive, and an explicit construction of $\CV(A_N)^{\bullet}$
will appear in Section~5.2.4 
below.

\pagebreak

\paragraph{5.2.~Proof.}  
\paragraph{5.2.1.~Beginning of the proof: a reduction to the homological degree
0.} 
Suppose one quantization, $\CV(A_N)^{\bullet}$, is given. We have the
direct sum decomposition
\[
\CV(A_N)^{\bullet}=\oplus_{n\in\BZ}\CV(A_N)^n,
\]
where $\CV(A_N)^0\subset \CV(A_N)$ is a vertex subalgebroid. Filtration (\ref{(2.8)}) in the
present situation becomes
\[
A_N\oplus\Omega(A_N)^0\subset\CV(A_N)^{\bullet}.
\]
Since $A_N\oplus\Omega(A_N)^0\subset\CV(A_N)^0$, $\CV(A_N)^0$ is a quantization of
$\CV^{\rm poiss}(A_N)^0$. This and the fact that $\Omega(A_N)^{\bullet}$ and $\Der(A_N)^{\bullet}$
are graded in opposite directions, cf.~(\ref{(3.1.8)}) and Lemma~4.5.1$(a)$, 
imply a~canonical vector space isomorphism
\begin{equation}
\CV(A_N)^{\bullet}\iso\CV(A_N)^{0}\bigoplus(\oplus_{n<0}\Der(A_N)^{n})\bigoplus(\oplus_{n>0}\Omega(A_N)^n).
\label{(4.2.1)}
\end{equation}
Now suppose that only $\CV(A_N)^{0}$ is given.

\medskip

\noindent
{\bf Lemma.} {\it  If $\CV$ is a quantization of $\CV^{\rm poiss}(A_N)^{0}$, then there is a unique
quantization, $\CV(A_N)^{\bullet}$, of $\CV^{\rm poiss}(A_N)^{\bullet}$ such that
$\CV(A_N)^0=\CV$.}

\begin{proof} $(i)$ {\it Uniqueness.} Pick a splitting (over $\BC$) of the exact sequence of graded vector
spaces, cf. (\ref{(2.9)}),
\[
0\rightarrow\Omega(A_N)^0\rightarrow\CV(A_N)_1^0\rightarrow \Der(A_N)^0\rightarrow 0
\]
so as to identify $\CV(A_N)_1$ with $\Omega(A_N)\oplus \Der(A_N)$ and obtain the projection
$\pi: \CV(A_N)_1\rightarrow\Omega(A_N)$ that is compatible with (\ref{(4.2.1)}). It follows
from Def\/inition~3.5 
that the only multiplications on $\CV(A_N)$ that are not
immediately determined by those on $\CV^{\rm poiss}(A_N)$ are the following components of ${}_{(0)}$
and ${}_{(1)}$:
\begin{gather*}
\pi\circ(_{(0)}):\ \  \Der(A_N)^{\bullet}\otimes
\Der(A_N)^{\bullet}\rightarrow\Omega(A_N)^{\bullet}, \qquad \xi\otimes\tau\mapsto\pi(\xi_{(0)}\tau),
\\
{}_{(1)}:\ \  \Der(A_N)^{\bullet}\otimes \Der(A_N)^{\bullet}\rightarrow A_N, \qquad
\xi\otimes\tau\mapsto\xi_{(1)}\tau.
\end{gather*}
The homological degree of the l.h.s.\ of these is non-positive, see (\ref{(4.2.1)}), of the
r.h.s.\ is non-negative; therefore, the operations may be non-zero only if both $\xi,\tau\in
\Der(A)^{0}$; the uniqueness follows.

$(ii)$ To prove the {\it existence}, note that
Def\/inition~3.2 
dif\/fers from Def\/inition~3.3 
in the following two respects only:

$\bullet$ the associativity of $_{(-1)}$ in the former is replaced with quasi-associativity
(\ref{(1.1e)}) in the latter;

$\bullet$ the requirements of the former that  multiplications ${}_{(n)}$, $n\geq 0$, be
derivations of  multiplication ${(-1)}$ and that multiplication ${}_{(-1)}$ be commutative are
simultaneously replaced with the Jacobi identity (\ref{(1.1d)}) with $m$ or $n$ equal to~$-1$.

Upon choosing a splitting as at the beginning of the proof, in each of the cases the identities
of Def\/inition~2.2 
exhibit {\it quantum corrections}, i.e., the terms that measure the failure
of a quantum object to be a classical one. It is easy to notice, by inspection, that in our
situation the quantum corrections may be non-zero only if all the terms involved belong to
$\CV=\CV(A_N)^0$.

One such example is provided by formula (\ref{(2.11b)}), where the failure of multiplication
${}_{(-1)}$ to be associative is measured by $-\tau(g)df-\tau(f)dg$; both the summands vanish
unless~$\tau\in\CV$, $f$,~$g$~being in $\CV$ automatically.

Another example deals with the commutativity of ${}_{(-1)}$. Let $f\in A_N$, $\tau\in
\Der(A_N)^{\bullet}$; then~(\ref{(1.1d)}) reads
\begin{gather*}
[\tau_{(-1)}, f_{(-1)}]=\tau(f)_{(-2)},
\end{gather*}
which is 0 unless $\tau\in\CV$. We leave it to the untiring reader to check the validity of all
the remaining requirements of Def\/inition~3.3. 
\end{proof}

In order to prove Theorem~5.1.1, 
it remains to quantize $\CV^{\rm poiss}(A_N)^0$. We
shall do this in a~somewhat roundabout manner.

\paragraph{5.2.2.~Localization.} 
To return to the hypothetical  vertex
algebroid $\CV(A_N)^{\bullet}$. By virtue of Section~5.2.1, 
it is enough to consider
$\CV(A_N)\stackrel{\text{def}}{=}\CV(A_N)^{0}$. Since $C={\rm Spec}(A_N)$
is af\/f\/ine, we can localize $\CV(A_N)$, see (\ref{(2.20)}) in
Section~3.7.1, 
so as to get a sheaf $\CV_{C}\in\Vert_{C}$.
Let $\check{C}$ be $C\setminus\{0\}$ and~$\CV_{\check{C}}$  the
restriction of~$\CV_{C}$ to~$\check{C}$. Apparently,
$\CV_{\check{C}}\in\Vert_{\check{C}}$ and, the manifold $\check{C}$
being smooth, our strategy will be to use the classif\/ication of the
objects of~$\Vert_{\check{C}}$, Section~3.7.3, 
so as to
identify those vertex algebroids over $\check{C}$ that may have come
from $C$ as above.

We begin by realizing  $\check{C}$  as a quotient of a manifold w.r.t.\ a f\/inite group action.
Consider the action
\begin{gather}
 \BZ_{N}\times\BC^{2}\rightarrow
\BC^{2},\qquad
 \bar{m}(y_1,y_2)=\big(\exp{(2\pi\sqrt{-1}m/N)}y_1,\exp{(2\pi\sqrt{-1}m/N)}y_2\big).
\label{(4.2.3)}
\end{gather}
The map
\begin{equation}
A_N\rightarrow\BC[y_1,y_2]^{\BZ_{N}},\qquad x_{j}\mapsto y_{1}^{N-j}y_{2}^{j}\label{(4.2.4a)}
\end{equation}
is an isomorphism; hence  isomorphisms
\begin{equation}
\BC^{2}/\BZ_{N}\iso C, \qquad (\BC^{2}\setminus 0)/\BZ_{N}\iso\check{C}. \label{(4.2.4b)}
\end{equation}
There arises  a projection
\begin{equation}
p: \ \ \BC^{2}\setminus 0\rightarrow\check{C}\label{(4.2.5)}
\end{equation}
and a faithful functor
\begin{equation}
p^*: \ \ \Vert_{\check{C}}\rightarrow \Vert_{\BC^{2}\setminus 0}\label{(4.2.6)}
\end{equation}
It is an equivalence of categories
\begin{equation}
p^*:\ \ \Vert_{\check{C}}\rightarrow \Vert_{\BC^{2}\setminus 0}^{\BZ_N},\label{(4.2.7)}
\end{equation}
where $\Vert_{\BC^{2}\setminus 0}^{\BZ_N}$ is the full subcategory of $\BZ_N$-equivariant
vertex algebroids; the inverse functor is that of $\BZ_N$-invariants: $\CV_{\BC^{2}\setminus
0}\mapsto \CV_{\BC^{2}\setminus 0}^{\BZ_N}$; cf.\ Section~3.7.1, 
(\ref{(2.21)}), (\ref{(2.22)}).

The objects of the category $\Vert_{\BC^{2}\setminus 0}$ were classif\/ied in
Section~3.7.3 
to the ef\/fect that there is a 1-1 correspondence between isomorphism
classes of vertex algebroids and linear combinations $\omega=\sum_{a,b>0}k_{ab}\omega_{ab}$,
where $\omega_{ab}$ is the 2-form $dy_1\wedge dy_2/y_1^ay_2^b$. It follows from the
construction that the vertex algebroid, $\CV_{\BC^{2}\setminus 0,\omega}$ is
$\BZ_N$-equivariant  if and only if $\omega$ is, hence if and only if
\begin{equation}
\omega=\sum_{a,b>0, N\text{ divides }a+b-2}k_{ab}\omega_{ab}.\label{(4.2.8)}
\end{equation}
It is from this list that we have to select.

\paragraph{5.2.3.~Conclusion of the proof.} 
By def\/inition, our hypothetical sheaf $\CV_{C}$ must  f\/it, for some~$\omega$, in
 the following commutative diagram:
\begin{equation}
\xymatrix{
0\ar[r]&\Gamma(C,\Omega_{C})\ar[r]\ar@{^{(}->}[d]&\Gamma(C,\CV_{C})\ar[r]\ar@{^{(}->}[d]&
\Gamma(C,\CT_{C})\ar[r]
\ar[d]&0\\
0\ar[r]&\Gamma(\check{C},\Omega_{\check{C}})\ar[r]&\Gamma(\check{C},\CV_{\check{C},\omega})\ar[r]&
\Gamma(\check{C},\CT_{\check{C}})\ar[r]&0
 } \label{(4.2.9)}
\end{equation}
Note that the vertical arrows are all the restriction (from $C$ to $\check{C}$) maps.
Furthermore, the rightmost vertical arrow is an equality. To see this, note that
$\Gamma(\check{C},\CT_{\check{C}})=\Gamma(\BC^2\setminus0,\CT_{\BC^2\setminus0})^{\BZ_N}$. The
latter is generated, over functions, by the tautological action of $gl_2$, see
Section~3.7.4, 
formu\-las~(\ref{(2.27)}),~(\ref{(2.28)}). (Indeed, an element of
$\Gamma(\BC^2\setminus0,\CT_{\BC^2\setminus0})^{\BZ_N}$ is a linear combination of
$f(y_1,y_2)\partial_1$ and $g(y_1,y_2)\partial_2$, where $N$ divides  $\deg(f)-1$ and
$\deg(g)-1$. This implies that $f(y_1,y_2)\partial_1$ is proportional to either $\rho(E_{11})$
or to $\rho(E_{21})$ and $g(y_1,y_2)\partial_2$ is proportional to either $\rho(E_{12})$ or to
$\rho(E_{22})$, see (\ref{(2.28)}).) Therefore, so is the former. But this action is precisely
the action of $gl_2$ on $C={\rm Spec}(A_N)$ described somewhat implicitly in (\ref{(4.1.1)}), hence
the equality $\Gamma(C,\CT_{ C})=\Gamma(\check{C},\CT_{\check{C}})$.

Contrary to this, the leftmost vertical arrow is not an equality; e.g. it is easy to check that
\begin{equation}
y_1^ry_2^{N-r-1}dy_1\in \Gamma(\check{C},\Omega_{\check{C}})\quad \text{but}\quad
y_1^ry_2^{N-r-1}dy_1\not\in \Gamma( C,\Omega_{C})\quad \text{if} \quad 0\leq r\leq N-2.\label{(4.2.10)}
\end{equation}
This simple remark is the reason why the quantization of $\CV^{\rm poiss}(A_N)$ is unique.

  Now, the upper row of (\ref{(4.2.9)}) is exact by virtue of
Def\/inition~3.5. 
This and the fact that the rightmost arrow is an equality imply that
the lower row must also be exact, at least on the right. By virtue of
Corollary~3.7.6, 
\begin{equation}
\omega=k\omega_{11}. \label{(4.2.11)}
\end{equation}
Now our task is to determine $k$.

Def\/ine $\CW$ to be the vertex subalgebroid of $\Gamma(\check{C},\CV_{\check{C},k\omega_{11}})$
generated by $A_N=\Gamma(\check{C},\CO_{\check{C}})$ and $\hat{\rho}(\CV(gl_2)_{-k-1,k-1})$,
where $\hat{\rho}$ is the one from Lemma~3.7.5. 

Since $\Gamma(C,\CT_{ C})$ is generated by $\rho(gl_2)$ over $A_N$,  $\CW=\Gamma(C,\CV_{C})$.
It is clear that $\CW=A_N\oplus\Gamma(C,\Omega_{C})+
A_{N(-1)}\hat{\rho}\CV^{\rm ch}(gl_2)_{-k-1,k-1}$ and were the elements $\rho(E_{ij})$ independent
over $A_N$, we would be done: $\CW$ would be the sought after quantization for any $k$. (In
fact, were that true, we could equivalently def\/ine $\CV(A_N)$ to be the push-out
$A_N\otimes\CV(gl_{2})_{-k-1,k-1}$, see (\ref{(2.23)}).) But they are not, and the  problem
with this is that an element of $A_{N(-1)}\hat{\rho}\CV^{\rm ch}(gl_2)_{-k-1,k-1}$ may belong to
$\Gamma(\check{C},\Omega_{\check{C}})$ and not to $\Gamma(C,\Omega_{C})$, cf.~(\ref{(4.2.10)}).

In fact, $A_N$ is a quadratic algebra, see (\ref{(4.1.0)}), and $\Gamma(C,\CT_{C})$ is a
quadratic $A_N$-module generated by $\{E_{ij},\;1\leq i,j\leq 2\}$. The relations, in terms of
$y_1$, $y_2$, are
\begin{equation}
y_{1}^{a}y_{2}^{b}y_{k}\hat{\rho}(E_{ij})-y_{1}^{a}y_{2}^{b}y_{i}\hat{\rho}(E_{kj})=0\quad \text{for all}\quad a+b+1=N,\label{(4.2.12)}
\end{equation}
as it easily follows from (\ref{(2.28)}).

Our task then is  to ensure that
\begin{equation}
(y_{1}^{a}y_{2}^{b}y_{k})_{(-1)}\hat{\rho}(E_{ij})-(y_{1}^{a}y_{2}^{b}y_{i})_{(-1)}\hat{\rho}(E_{kj})
\in\Gamma(C,\Omega_{C}) \quad \text{for all} \quad a+b+1=N.\label{(4.2.13)}
\end{equation}

Let us consider for the sake of def\/initeness the case of $i=1$, $k=j=2$. We have to compute
 the following section of $\Gamma(\check{C},\CV_{\check{C},k\omega_{11}}(\check{C}))$:
\begin{equation}
(y_{1}^{r}y_{2}^{N-r})_{(-1)}\hat{\rho}(E_{12})-(y_{1}^{r+1}y_{2}^{N-r-1})_{(-1)}\hat{\rho}(E_{22})\quad \text{for all} \quad  0\leq r<N.\label{(4.2.14)}
\end{equation}
Formulas (\ref{(2.29a)})--(\ref{(2.29d)}) allow us to re-write
(\ref{(4.2.14)}) as follows
\[
(y_{1}^{r}y_{2}^{N-r})_{(-1)}(y_{1}\partial_{2})-(y_{1}^{r+1}y_{2}^{N-r-1})_{(-1)}
\left(y_2\partial_{y_{2}}+\frac{kT(y_{1})}{y_1}\right).
\]
A little thought (or formula (\ref{(2.11b)})) will show that the f\/irst term yields
\begin{gather*}
 y_{1}^{r+1}y_{2}^{N-r}\partial_{2}
- r(N-r)y_{1}^{r}y_{2}^{N-r-1}T(y_{1}) -(N-r)(N-r-1)y_{1}^{r+1}y_{2}^{N-r-2}T(y_{2}) ,
\end{gather*}
where it is understood that
\[
y_{1}^{r+1}y_{2}^{N-r}\partial_{2}\stackrel{\text{def}}{=}(y_{1(-1)}(y_{1(-1)}(\cdots
(y_{1(-1)}(y_{2(-1)}(\cdots(y_{2(-1)}\partial_{2})\cdots)))))).
\]

 The second one will likewise give
\begin{gather*}
 -y_{1}^{r+1}y_{2}^{N-r}\partial_{2}
 +(N-r-1)(r+1)y_{1}^{r}y_{2}^{N-r-1}T(y_{1})\\
 \qquad{}
+(N-r-1)(N-r-2)y_{1}^{r+1}y_{2}^{N-r-2}T(y_{2})
 -ky_{1}^{r}y_{2}^{N-r-1}T(y_{1}) .
\end{gather*}

Adding one to another makes  expression (\ref{(4.2.14)}) into
\[
(N-1-2r-k)y_{1}^{r}y_{2}^{N-r-1}T(y_{1}) +(-2N+2r+2)y_{1}^{r+1}y_{2}^{N-r-2}T(y_{2}).
\]
The latter equals the total derivative
\[
-2T(y_{1}^{r+1}y_{2}^{N-r-1})=-2T(x_{N-r-1}),
\]
and is therefore an element of $\Gamma(C,\Omega_{C})$, precisely when $k=N+1$.

The case where $i=j=1$, $k=2$ works out similarly and gives the same answer $k=N+1$. This
concludes the proof of item (a) of Theorem~5.1.1. 

As to item $(b)$, (\ref{(4.1.4)}) follows from the $k=N+1$ case of Lemma~3.7.5, 
the
assertion that has been instrumental for the proof anyway, and (\ref{(4.1.5)}) follows from
(\ref{(2.32)}). \hfill $\qed$

\paragraph{Corollary 5.2.4.} 
{\it The unique quantization
$\CV(A_N)=\CV(A_N)^0\oplus \Der(A_N)^{\bullet}\oplus\Omega(A_N)^{\bullet}$, where $\CV(A_N)^0$
is the vertex subalgebroid of $\Gamma(\check{C},\CV_{\check{C},(N+1)\omega_{11}})$ generated by
$A_N=\Gamma(\check{C},\CO_{\check{C}})$ and \linebreak $\hat{\rho}(\CV(gl_2)_{-N-2,N})$, $\hat{\rho}$
being the one from Lemma~{\rm 3.7.5.}} 

\paragraph{5.3.~Higher dimensional Veronese embeddings.} 
Regard $\BC^n$ as the tautological representation of $gl_n$ and let $V=(\BC^n)^*$. Let
\begin{equation}
\iota_N: \ \ \BP(V)\rightarrow \BP(S^N V), \quad l\mapsto l^{\otimes N} \label{(4.3.1)}
\end{equation}
be the classical Veronese embedding. By $A_{nN}$ let us denote the homogeneous coordinate ring
of $\iota_N(\BP(V))$. It is clear that if $n=2$, then $A_{nN}$ is the algebra $A_N$ we  dealt
with above. It is now natural to ask if $\CV^{\rm poiss}(A_{nN})$ af\/fords a quantization. The
result is a bit disheartening.

\paragraph{Theorem 5.3.1.} 
{\it The vertex Poisson algebroid
$\CV^{\rm poiss}(A_{nN})$ cannot be quantized if $N>1$ and $n>2$.}

\begin{proof} Consider the action
\begin{equation}
\BZ_N\times V\rightarrow V,\qquad (\bar{m},v)\mapsto \exp{(2\pi\sqrt{-1}m/N)}v.\label{(4.3.2)}
\end{equation}
Analogously to (\ref{(4.2.4a)}), (\ref{(4.2.4b)}), we obtain isomorphisms
\begin{equation}
V/\BZ_N\iso {\rm Spec}(A_{nN}),\qquad (V\setminus 0)/\BZ_N\iso {\rm Spec}(A_{nN})\setminus 0.\label{(4.3.3)}
\end{equation}
Thus we are led to the question, ``How many vertex algebroids are there on  $(V\setminus
0)/\BZ_N$?'' That such vertex algebroids exist is obvious because $\CV_{V\setminus 0}^{\BZ_N}$
is one; here $\CV_{V\setminus 0}$ is the pull-back of the standard $\CV_{V}$, cf.\
Section~3.7.2, 
on $V\setminus 0$. Note that if $x_1,\dots ,x_n$ is a basis of $\BC^n$~--
remember that we think of $\BC^n$ as the space of linear functions on~$V$~-- then the
assignment
\begin{equation}
\hat{\rho}: \ \ E_{ij}\mapsto x_{i(-1)}\partial_j \label{(4.3.4a)}
\end{equation}
def\/ines a vertex algebroid morphism
\begin{equation}
\hat{\rho}: \ \ \CV(gl_n)_{-1,-1}\rightarrow \Gamma\big((V\setminus 0)/\BZ_N,\CV_{V\setminus
0}^{\BZ_N}\big).\label{(4.3.4b)}
\end{equation}

\paragraph{Lemma 5.3.2.} 
{\it The manifolds $V\setminus 0$,
$(V\setminus 0)/\BZ_N$ carry a unique up to isomorphism sheaf of vertex algebroids. It is
isomorphic to $\CV_{V\setminus 0}$ in the former case and to $\CV_{V\setminus 0}^{\BZ_N}$ in
the latter case.}

\begin{proof}[Proof of Lemma~5.3.2.] 
By virtue of
 (\ref{(2.17)}) and (\ref{(2.21)}), it suf\/f\/ices to show that
\begin{equation}
H^1\big(V\setminus 0,\Omega^2_{V\setminus 0}\rightarrow \Omega^{3,{\rm cl}}_{V\setminus
0}\big)=0.\label{(4.3.5)}
\end{equation}
Converging to the hypercohomology $H^{*}(V\setminus 0,\Omega^2_{V\setminus 0}\rightarrow
\Omega^{3,{\rm cl}}_{V\setminus 0})$  is a standard spectral sequence with
\begin{gather*}
E^1_{00}=H^0\big(V\setminus 0,\Omega^2_{V\setminus 0}\big),
\\
E^1_{01}\oplus E^1_{10}=H^1\big(V\setminus 0,\Omega^2_{V\setminus 0}\big)\oplus H^0\big(V\setminus 0,
\Omega^{3,{\rm cl}}_{V\setminus 0}\big).
\end{gather*}
The next dif\/ferential is
\[
d^2=d_{DR}: \ \ H^0\big(V\setminus 0,\Omega^2_{V\setminus 0}\big)\rightarrow H^0\big(V\setminus 0,
\Omega^{3,{\rm cl}}_{V\setminus 0}\big),
\]
and it is clear that it is surjective. Finally, if $n>2$, then $H^1(V\setminus
0,\Omega^2_{V\setminus 0})=0$; this concludes the proof of Lemma~5.3.2. 
\end{proof}

It is clear now why $\CV^{\rm poiss}(A_{nN})$ cannot be quantized if $n>2$: it is because
Lemma~5.3.2 
has left us no room for manoeuvre that was helpfully provided by the
analysis of Section~3.7.3 
in the $n=2$ case. Indeed, one can now repeat the entire
argument of Sections~5.2.1--5.2.3 
only to f\/ind out that an obvious analogue
of (\ref{(4.2.13)}) is false. Here are some details:

According to Lemma~5.3.2, 
a quantization of $\CV^{\rm poiss}(A_{nN})$, if existed, upon
localizing to $(V\setminus 0)/\BZ_N$ would give $\CV_{V\setminus 0}^{\BZ_N}$. Hence any such
quantization can  be equal only to the vertex subalgebroid of $\Gamma(V\setminus
0,\CV_{V\setminus 0})^{\BZ_N}$ generated by $A_{nN}$ and $\hat{\rho}(\CV(gl_n)_{-1,-1}$, see
(\ref{(4.3.4a)}-\ref{(4.3.4b)}). But this subalgebroid necessarily contains elements from
$\Gamma((V\setminus 0)/\BZ_N,\Omega_{(V\setminus 0)/\BZ_N})\setminus \Omega(A_{nN})$. Indeed, a
computation analogous to the one performed at the end of Section~5.2.3 
shows that if
$n\geq 3$, then
\begin{gather*}
 \big(x_3x_2^{N-1}\big)_{(-1)}(x_{1(-1)}\partial_{1})-\big(x_3x_2^{N-2}x_1\big)_{(-1)}(x_{2(-1)}\partial_{1})\\
 \qquad {} =\big(T\big(x_3x_2^{N-2}\big)\big)_{(-1)}x_2
 \in \Gamma((V\setminus 0)/\BZ_N,\Omega_{(V\setminus 0)/\BZ_N}))\setminus
\Omega(A_{nN}) .
\end{gather*}
This concludes the proof of Theorem~5.3.1. 
\end{proof}

\section{Chiral Hamiltonian reduction
interpretation}\label{sect.5}

 We will now interpret some of the
constructions above in the language of semi-inf\/inite cohomology. Our
exposition will be brief and almost no proofs will be given. In some
respects, the material of this section is but an afterword to~\cite{GMSII}.

Since we will be mostly concerned with smooth varieties, we will
f\/ind it convenient to work not with vertex algebroids, such as
$\CV(\fg)_k$, $\CV_{X}$, but with the corresponding vertex algebras
or CDO-s, such as $\CU_{\rm alg}\CV(\fg)_k=\CU\CL(\fg)_k$,
$D^{\rm ch}_{X}=\CU_{\rm alg}\CV_{X}$.

\paragraph{6.1.~Semi-inf\/inite cohomology.} 
Let
$\BV$ be a vertex algebra, $\fg$ a f\/inite dimensional Lie algebra,
$(\cdot,\cdot)$ an invariant bilinear form on $\fg$, and $\rho$ a vertex
algebra morphism
\begin{equation}
\rho:\ \  \CU_{\rm alg}\CV(\fg)_{(\cdot,\cdot)}\rightarrow\BV,\label{(5.1.1)}
\end{equation}
see Example~2.5,  especially (\ref{(1.5c)}).

Introduce the Clif\/ford vertex algebra built on
$\Pi(\fg\oplus\fg^*)$, $\Pi$ being the parity change functor. This
vertex algebra is nothing but the space of global sections of the
standard CDO on superspace $\Pi(\fg\oplus\fg^*)$, see
Section~3.7.2 
for a discussion of a purely even analogue.
Denote this vertex algebra by $D^{\rm ch}(\Pi(\fg\oplus\fg^*))$.

By def\/inition, if we let $\{x_i\}$ be a basis of $\fg$, $\{\phi_i\}$
the corresponding basis of $\Pi(\fg)$,  $\{\phi_i^*\}$ the dual
basis of $\Pi(\fg^*)$, then $D^{\rm ch}(\Pi(\fg\oplus\fg^*))$ is the
vertex algebra  generated by the vector space $\Pi(\fg\oplus\fg^*)$
and relations
\begin{gather*}
(\phi_i^*)_{(0)}\phi_j=(\phi_i)_{(0)}\phi_j^*=\delta_{ij}\b1,
\\
(\phi_i^*)_{(n+1)}\phi_j=(\phi_i)_{(n+1)}\phi_j^*=(\phi_i^*)_{(n)}\phi^*_j=(\phi_i)_{(n)}\phi_j^*=0
\quad \text{if} \quad n\geq 0.
\end{gather*}
There arises the vertex algebra $\BV\otimes
D^{\rm ch}(\Pi(\fg\oplus\fg^*))$.

If $\{c_{ij}^k\}$ are the structure constants of $\fg$ relative to
$\{x_i\}$, that is to say, if $[x_i,x_j]=\sum_kc_{ij}^{k}x_k$, then
following~\cite{Feig} one considers the element
\[
d^{\infty/2}=\sum_k\rho(x_k)_{(-1)}\phi^*_k-
\frac{1}{2}\sum_{i,j,k}c_{ij}^{k}\phi_{k(-1)}(\phi^*_{i(-1)}\phi^*_j)\in
\BV\otimes D^{\rm ch}(\Pi(\fg\oplus\fg^*)).
\]
A direct computation shows that
\begin{equation}
(d^{\infty/2})_{(0)}d^{\infty/2}=0 \quad \text{if} \quad
(\cdot,\cdot)=-K(\cdot,\cdot),\label{(5.1.2)}
\end{equation}
where $K(\cdot,\cdot)$ is the Killing form on $\fg$: $K(a,b)={\rm tr} ({\rm ad}_a\cdot
{\rm ad}_b)$.

If condition (\ref{(5.1.2)}) is satisf\/ied, then we obtain a DGVA
$C^{\infty/2}(L\fg;\BV)\stackrel{\text{def}}{=} (\BV\otimes
D^{\rm ch}(\Pi(\fg\oplus\fg^*))$, $d^{\infty/2}_{(0)})$. The cohomology
vertex algebra $H^{\infty/2}(L\fg;\BV)$ is due to Feigin \cite{Feig}
and well known as either semi-inf\/inite or BRST cohomology of the
loop algebra $L\fg$ with coef\/f\/icients in $\BV$. If one chooses to
think of $\BV$ as an algebra of (`chiral') functions on a symplectic
manifold with $\fg$-structure, then $H^{\infty/2}(L\fg,\BV)$ is to
be thought of as an algebra of functions on the symplectic quotient
$M//\fg$, hence the title of this section.

One similarly def\/ines the relative version
$H^{\infty/2}(L\fg,\fg;\BV)$, see \cite{GMSII} for some details; the
condition (\ref{(5.1.2)}) remains the same in this case.

\paragraph{6.2.~The $\boldsymbol{sl_2}$ case.} 
Let us return to the set-up of Section~5.2.2, 
where we had
the Veronese cone $C={\rm Spec}(A_N)$, $\check{C}=C\setminus 0$, and
consider $\CL_N$, the degree $N$ line bundle over $\BP^1$, and
$\check{\CL}_N=\CL_N\setminus\{\text{the zero section}\}$. We
obtain the commutative square
\begin{equation}
\xymatrix{ \CL_N\ar[r]&C\\
  \check{\CL}_N\ar[r]^{\sim}\ar@{^{(}->}[u]&\check{C}\ar@{^{(}->}[u]
 }
\label{(5.2.1)}
\end{equation}
where the upper horizontal map is a surjective birational
isomorphism, a blow-up of the vertex of the cone. We have seen  that
$\check{C}$ carries a family of CDO-s, $D^{\rm ch}_{\check{C},\omega}$,
$\omega\in H^1(\check{C},\Omega^2_{\check{C}})$. Denote the
coorresponding family of CDO-s on $\check{\CL}_N$ by
$D^{\rm ch}_{\check{\CL}_N,\omega}$. Theorem~5.1.1 
says that
$D^{\rm ch}_{\check{C},\omega}$ is a pull-back of a CDO on $C$ if\/f
$\omega=(N-1)\omega_{11}$, in which case it contains
$\CV(sl_2)_{-N-2}$.

Now a question arises, ``For what, if any, $\omega$ is
$D^{\rm ch}_{\check{\CL}_N,\omega}$ a pull-back of a CDO from $\CL_N$?''
The existence of such $\omega$ depends on the vanishing of the
characteristic class ${\rm ch}_2(\CL_N)$, (\ref{(2.16)}). A simple way to
prove the vanishing result, and to compute a possible $\omega$, is
provided by the semi-inf\/inite cohomology.

One has  the base af\/f\/ine space $\BC^2\setminus 0$, the principal
$N$-bundle $p:\;SL_2\rightarrow \BC^2\setminus 0$, where $N$ is the
subgroup of upper-triangular matrices and a family of CDO-s,
$\CD_{SL_2,(\cdot,\cdot)}$ on $SL_2$, over $SL_2$. This family enjoys
\cite{AG,GMSII} the 2 vertex algebra embeddings
\begin{gather}
\CU_{\rm alg}\CV(\fg)_{(\cdot,\cdot)}\stackrel{\rho_l}{\rightarrow}
\Gamma(SL_2,D^{\rm ch}_{SL_2,(\cdot,\cdot)})\stackrel{\rho_r}{\leftarrow}
\CU_{\rm alg}\CV(\fg)_{-(\cdot,\cdot)-K(\cdot,\cdot)},\nonumber\\
\text{s.t.} \quad \rho_l(a)_{(n)}\rho_r(b)=0\quad \text{if} \quad n\geq 0. \label{(5.2.2)}
\end{gather}

  In this case condition (\ref{(5.1.2)})
is  satisf\/ied for all forms $(\cdot,\cdot)$. Therefore, for any $U\subset
SL_2/N$, there arises a vertex algebra $H^{\infty}(L\fn,
\Gamma(p^{-1}(U),D^{\rm ch}_{SL_2,k}))$, where we use $\rho_r$ in place
of  $\rho$, see (\ref{(5.1.1)}). Denote by $H^{\infty/2}(L\fn,
D^{\rm ch}_{SL_2,(\cdot,\cdot)})$ the sheaf associated with the presheaf
$U\mapsto H^{\infty}(L\fn, \Gamma(p^{-1}(U)$, $\CD_{SL_2,(\cdot,\cdot)}))$. It
was noted in \cite{GMSII} that this sheaf is a CDO, which re-proves
the obvious fact that ${\rm ch}_2(\BC^2\setminus 0)=0$.

Note that the left one of the embeddings (\ref{(5.2.2)}) survives
the passage to the cohomology. The right one does not, not entirely
at least, but the embedding of the torus part does, albeit with a
shifted central charge. We obtain
\begin{equation}
\CU_{\rm alg}\CV(\fg)_{(\cdot,\cdot)}\stackrel{\rho_l}{\rightarrow}
\Gamma\big(\BC^2\setminus 0,H^{\infty/2}(L\fn,
D^{\rm ch}_{SL_2,(\cdot,\cdot)})\big)\stackrel{\rho_r}{\leftarrow}
\CU_{\rm alg}\CV(\ft)_{-(\cdot,\cdot)|_{\ft}-1/2K(\cdot,\cdot)|_{\ft}} \label{(5.2.3)}
\end{equation}
so that $ \rho_l(a)_{(n)}\rho_l(b)=0\text{ if }n\geq 0$; here
$(\cdot,\cdot)|_{\ft}$ and $K(\cdot,\cdot)|_{\ft}$ stand for the restrictions of the
corresponding forms to $\ft$.

Altogether, the 2 embeddings provide a vertex algebra embedding of
$\CU_{\rm alg}\CV(gl_2)_{k_1,k_2}$ with appropriate central charges
$k_1$ and $k_2$. In fact, if we let $(a,b)=k\,{\rm tr}(a\cdot b)$ as we  did
in  Section~3.7.4, 
then we obtain that $k_1=k$, $k_2=-k-2$
and a diagram, cf. Lemma~3.7.5 
and (\ref{(2.30)}).
\begin{equation}
\CU_{\rm alg}\CV(\fg)_k\hookrightarrow \Gamma(\BC^2\setminus
0,H^{\infty/2}(L\fn, D^{\rm ch}_{SL_2,(\cdot,\cdot)}))\iso \Gamma(\BC^2\setminus
0, D^{\rm ch}_{\BC^2\setminus 0,-(k+1)\omega_{11}}).\label{(5.2.4)}
\end{equation}
Lemma~3.7.5 
shows that the chiral hamiltonian reduction
technology reproduces precisely those CDO-s on the punctured plane
that carry an af\/f\/ine Lie algebra action.

In order to try and def\/ine a CDO on $\CL_N$, we represent the latter
as
\begin{equation}
\CL_N=(\BC^2\setminus 0)\times_{\BC^*}\BC_N,\label{(5.2.5)}
\end{equation}
where $\BC_N$ is the character $\BC^*\ni z \mapsto z^N$. This
suggests to def\/ine a CDO on $\CL_N$ as the chiral hamiltonian
reduction of $D^{\rm ch}_{\BC^2\setminus 0,-(k+1)\omega_{11}}\otimes
D^{\rm ch}_{\BC}$ using embedding (\ref{(5.2.3)}) twisted by action
(\ref{(5.2.5)}); in practice that means that if $h$ is the standard
generator of $\ft$, then one has to replace $\rho_r(h)$ in
(\ref{(5.2.3)}) with $\rho_r(h)+ Ny_{(-1)}\partial_y$, where $y$ is
a coordinate on $\BC_N$. Two things are to be kept in mind: f\/irst,
since the topology of $\BC^*$ is non-trivial \cite{GMSII}, one has
to use the relative version of the semi-inf\/inite cohomology; second,
and most important, condition (\ref{(5.1.2)}) is not automatically
satisf\/ied. In fact, (\ref{(5.1.2)}) is equivalent to
\begin{equation}
(\rho_r(h)+ Ny_{(-1)}\partial_y)_{(1)}(\rho_r(h)+
Ny_{(-1)}\partial_y)=0,\label{(5.2.6a)}
\end{equation}
which gives
\begin{equation}
(\cdot,\cdot)=-\left(\frac{N^2}{8}+\frac{1}{2}\right)K(\cdot,\cdot).\label{(5.2.6b)}
\end{equation}
It follows that $\CH^{\infty/2}(L\ft,\ft;D^{\rm ch}_{\BC^2\setminus
0,-(k+1)\omega_{11}}\otimes D^{\rm ch}_{\BC})$ is well def\/ined and gives
a CDO on $\CL_N$ precisely if (\ref{(5.2.6b)}) holds.

To conclude,
\begin{enumerate}\itemsep=0pt
\item[(1)] the manifold $\check{\CL}_N$ carries a 1-parameter family of
CDO-s, $D^{\rm ch}_{\BC^2\setminus 0,-(k+1)\omega_{11}}$, with \linebreak
$\CU_{\rm alg}\CV(sl_2)_{(\cdot,\cdot)}$-structure, see Section~3.7.4; 

\item[(2)] the condition that a CDO on $\check{C}$ extends to one on $C$
picks  a unique
  $(\cdot,\cdot)$; the latter depends on $N$
linearly, see Theorem~5.1.1; 

\item[(3)] the family $D^{\rm ch}_{\BC^2\setminus 0,-(k+1)\omega_{11}}$
contains at least one representative that extends to a CDO on~$\CL_N$; the condition that a CDO on $\check{\CL}_N$ extends to
$\CL_N$ and af\/fords a realization via the chiral Hamiltonian
reduction picks  a unique  $(\cdot,\cdot)$; the latter depends on $N$
quadratically,~(\ref{(5.2.6b)}).
\end{enumerate}

In fact, there is a third way to f\/ix a $(\cdot,\cdot)$. This one amounts to
carrying a regularization procedure \`a la Lambert, used in~\cite{BN} in
a similar but dif\/ferent situation, and gives another quadratic
dependence on $N$. The importance of this approach is yet to be
worked out.

\paragraph{6.3.~Higher rank generalization.} 
Let $G$ be a simple complex Lie group, $P\subset G$ a para\-bo\-lic
subgroup, $R\subset P$ the unipotent radical of $R$, $M=P/R$. Let us
make the following assumption
\begin{equation}
M=M_0\times M_1,\qquad \text{where} \quad M_0\iso\BC^*,\quad  M_1 \ \ \text{is
simple.}\label{(5.3.1)}
\end{equation}
Let $Q\subset G$ be the extension of $M_1$ by $R$. Thus we obtain a
$\BC^*$-bundle
\[
G/Q\rightarrow G/P
\]
and the associated line bundle
\[
\CL_Q\rightarrow G/P.
\]
We have
\begin{equation}
{\rm ch}_2(G/Q)=0.\label{(5.3.2)}
\end{equation}
Indeed, there is at least one CDO on $G/Q$ that can be def\/ined via
the chiral Hamiltonian reduction as follows. By analogy with
Section~6.2, 
since (\ref{(5.2.2)}) holds true with $SL_2$
replaced with an arbitrary simple complex Lie group $G$, we observe
that if $\fq=Lie(Q)$, then $H^{\infty/2}(L\fq,\fq;D^{\rm ch}_{G,(\cdot,\cdot)})$
is well def\/ined for precisely one choice of $(\cdot,\cdot)$. In fact, in
this case condition (\ref{(5.1.2)}) amounts to the requirement that
the restriction $(\cdot,\cdot)$ to $M_1$ be equal to the Killing form on
$M_1$, and there is one and only one way to achieve that by
appropriately rescaling
 $(\cdot,\cdot)$  -- this
is where assump\-tion~(\ref{(5.3.1)}) is crucial. It is rather clear
that $H^{\infty/2}(L\fq,\fq;D^{\rm ch}_{G,(\cdot,\cdot)})$ is a CDO on $G/Q$.

As a consequence, we obtain a vertex algebra morphism
\begin{equation}
\CU_{\rm alg}\CV(\fg)_{(\cdot,\cdot)}\rightarrow
H^{\infty/2}(L\fq,\fq;D^{\rm ch}_{G,(\cdot,\cdot)})\label{(5.3.3)}
\end{equation}
for a uniquely determined bilinear form $(\cdot,\cdot)$.

Assertion (\ref{(5.3.2)}) is of course analogous to the fact that
${\rm ch}_2(SL_2/N)=0$, which was discussed in Section~6.2. 
Unlike the $sl_2$-case, however, we have  obtained not  a family but
a single CDO on $G/Q$, and this precludes the def\/inition of a CDO on
$\CL_Q$ via the chiral Hamiltonian reduction. It is then natural to
expect that  ${\rm ch}_2(\CL_Q)\neq 0$ and that each CDO on
 $G/Q$ with $\CU_{\rm alg}\CV(\fg)_{(\cdot,\cdot)}$-structure is isomorphic to
  $H^{\infty/2}(L\fq,\fq;D^{\rm ch}_{G,(\cdot,\cdot)})$.

One example of this analysis is provided by $G=SL_n$ with $Q$ chosen
to be the subgroup with 1st column equal to $(1,0,0,\dots ,0)$. Then
$G/Q$ is $\BC^n\setminus 0$, the corresponding CDO has been used in
Section~5.3, 
and it is easy to check that in this case the
embedding (\ref{(5.3.3)}) becomes precisely
(\ref{(4.3.4a)}), (\ref{(4.3.4b)}).

Another example is provided by the {\it space of pure spinors}
punctured at a point.
 It is a~homogeneous space which satisf\/ies assumption (\ref{(5.3.1)}).
  Therefore, our analysis is an alternative
 way to prove the vanishing of the 2nd component
of the Chern character, originally verif\/ied by Nekrasov~\cite{N}.

Needless to say, our discussion is very close  in spirit to the
def\/inition of Wakimoto modules due to Wakimoto and Feigin--Frenkel,
see \cite{F2} and references therein. In fact, it is easy to see
that the spaces of sections over `the big cell' of the sheaves
constructed in Section~6.2 
contain Wakimoto modules over
$\widehat{sl}_2$, and those of the present section contain the
so-called {\it generalized Wakimoto modules corresponding to}~$R$.

\subsection*{Acknowledgements}

The author would like to thank N.~Nekrasov and especially V.~Hinich
for interesting discussions and for bringing~\cite{Hin} to his attention.
The paper was completed at the IHES in Bures-sur-Yvette. We are
grateful to the institute for hospitality and excellent working
conditions. This work  was partially supported by an NSF grant.
Special thanks go to V.~Gorbounov and V.~Schechtman.

\pdfbookmark[1]{References}{ref}
\LastPageEnding


\begin{thebibliography}{99}

\footnotesize\itemsep=0pt

\bibitem{AG} Arkhipov S., Gaitsgory D., Dif\/ferential operators on the loop
group via chiral algebras, {\it  Int. Math. Res. Not.}  {\bf 2002} (2002), no.~4,
165--210, \href{http://arxiv.org/abs/math.AG/0009007}{math.AG/0009007}.

\bibitem{Behr} Behrend K., Dif\/ferential graded schemes I: Perfect resolving
algebras, \href{http://arxiv.org/abs/math.AG/0212225}{math.AG/0212225}.

\bibitem{BF} Backelin J., Fr\"oberg R., Koszul algebras, Veronese subrings
and rings with linear resolutions, {\it Rev. Roumaine Math. Pures
Appl.} {\bf 30} (1985), no.~2, 85--97.

\bibitem{BN} Berkovits N., Nekrasov N., The character of pure spinors, {\it Lett. Math. Phys.} {\bf
74} (2005), 75--109, \mbox{\href{http://arxiv.org/abs/hep-th/0503075}{hep-th/0503075}}.

\bibitem{Ber} Berkovits N., Super-Poincar\'e covariant quantization of the
superstring, {\it J. High Energy Phys.} {\bf 2000} (2000), no.~4, 18, 17~pages, \href{http://arxiv.org/abs/hep-th/0001035}{hep-th/0001035}.

\bibitem{Bez} Bezrukavnikov R., Koszul property and Frobenius splitting of
Schubert varieties, \href{http://arxiv.org/abs/alg-geom/9502021}{alg-geom/9502021}.

\bibitem{Bre} Bressler P., The f\/irst Pontryagin class, {\it Compos. Math.} {\bf 143} (2007),
1127--1163, \href{http://arxiv.org/abs/math.AT/0509563}{math.AT/0509563}.

\bibitem{C} Courant T.J., Dirac manifolds, {\it Trans. Amer. Math. Soc.} {\bf
319} (1990), 631--661.

\bibitem{Dor} Dorfman I., Dirac structures of integrable evolution
equations, {\it Phys. Lett.~A} {\bf 125} (1987), no.~5, 240--246.

\bibitem{F1} Frenkel~E., Private communication.

\bibitem{F2} Frenkel~E., Wakimoto modules, opers and the center at the critical level, {\it  Adv. Math.} {\bf  195}  (2005),  297--404,
\href{http://arxiv.org/abs/math.QA/0210029}{math.QA/0210029}.

\bibitem{FBZ} Frenkel E., Ben-Zvi D.,  Vertex algebras and algebraic curves, 2nd ed.,
{\it Mathematical Surveys and Monographs}, Vol.~88, American Mathematical Society, Providence, RI, 2004.

\bibitem{Feig} Feigin B., Semi-inf\/inite homology of Lie, Kac--Moody and
Virasoro algebras, {\it Uspekhi Mat. Nauk} {\bf  39} (1984),  no.~2, 195--196 (in Russian).

\bibitem{FP} Feigin B., Parkhomenko S., Regular representation of af\/f\/ine
Kac--Moody algebras, in Algebraic and geometric methods in
mathematical physics (Kaciveli, 1993), {\it Math. Phys. Stud.},
Vol.~19, Kluwer Acad. Publ., Dordrecht, 1996, 415--424, \href{http://arxiv.org/abs/hep-th/9308065}{hep-th/9308065}.

\bibitem{GMSI} Gorbounov V., Malikov F., Schechtman V.,  Gerbes of chiral
dif\/ferential operators. II.~Vertex algebroids, {\em Inv. Math.} {\bf
155} (2004), 605--680, \href{http://arxiv.org/abs/math.AG/0003170}{math.AG/0003170}.

\bibitem{GMSII} Gorbounov V., Malikov F., Schechtman V.,  On chiral
dif\/ferential operators over homogeneous spaces, {\em Int. J. Math. Math. Sci.} {\bf 26}
(2001), no. 2, 83--106, \href{http://arxiv.org/abs/math.AG/0008154}{math.AG/0008154}.

\bibitem{GS} Gorbounov V., Schechtman V.,
 Homological algebra and divergent series,
\href{http://arxiv.org/abs/arXiv:0712.3670}{arXiv:0712.3670}.




\bibitem{Hin} Hinich V., Homological algebra of homotopy algebras,
 {\it Comm. Algebra}  {\bf 25}  (1997),  3291--3323, \href{http://arxiv.org/abs/q-alg/9702015}{q-alg/9702015}.

\bibitem{IM} Inamdar S.P., Mehta V.B., Frobenius splitting of Schubert
varieties and linear syzygies, {\it Amer. J. Math.} {\bf 116}
(1994), 1569--1586.

\bibitem{Kac} Kac V., Vertex algebras for beginners, 2nd ed., {\it University Lecture Series}, Vol.~10, American Mathematical Society, Providence, RI, 1998.


\bibitem{Li} Li H., Abelianizing vertex algebras,
 {\it Comm. Math. Phys.} {\bf  259}  (2005),   391--411,
\href{http://arxiv.org/abs/math.QA/0409140}{math.QA/0409140}.

\bibitem{LWX} Liu Z.-J., Weinstein A., Xu P., Manin triples for Lie
bialgebroids, {\it J. Differential Geom.} {\bf 45} (1997), 547--574, \href{http://arxiv.org/abs/dg-ga/9508013}{dg-ga/9508013}.

\bibitem{Mal} Malikov F., Lagrangian approach to sheaves of vertex algebras,
{\it Comm. Math. Phys.}  {\bf 278} (2008), 487--548, \href{http://arxiv.org/abs/math.AG/0604093}{math.AG/0604093}.


\bibitem{N} Nekrasov N., Lectures on curved beta-gamma systems, pure spinors, and anomalies,
\href{http://arxiv.org/abs/hep-th/0511008}{hep-th/0511008}.

\bibitem{Pr} Primc M., Vertex algebras generated by Lie algebras, {\it J.
 Pure Appl. Algebra} {\bf 135} (1999), 253--293, \href{http://arxiv.org/abs/math.QA/9901095}{math.QA/9901095}.

\end{thebibliography}
\end{document}